\documentstyle[amssymb,graphicx,12pt]{article}

\input tcilatex

\begin{document}

\title{Modelling of thin isotropic elastic plates with small
piezoelectric inclusions and distributed electric circuits. Models
for inclusions larger or comparable to the thickness of the plate.}

\author{ \'{E}ric Canon \\
{\small Pres universit\'{e} de Lyon,}\\
[-0.8ex] {\small UMR CNRS 5208, }\\
[-0.8ex] {\small Universit\'{e} Jean Monnet de Saint-\'{E}tienne, }\\
{\small Institut Camille Jordan,}\\
[-0.8ex] {\small 23\ rue du docteur Paul Michelon, F-42\ 023\ Saint-%
\'{E}tienne. France.}\\
[-0.8ex] {\small \texttt{email: eric.canon@univ-st-etienne.fr}}\\
\and Michel Lenczner \\
{\small FEMTO-ST, D\'{e}partement Temps-Fr\'{e}quence,}\\
[-0.8ex] {\small University of Franche-Comt\'{e},}\\
[-0.8ex] {\small 26 chemin de l'Epitaphe, 25030 Besan\c{c}on Cedex, France}\\
[-0.8ex] {\small \texttt{michel.lenczner@utbm.fr} }}

\maketitle

\bigskip

{\bf \bigskip \noindent Abstract. \ }This paper is the second part of a work
devoted to the modelling of thin elastic plates with small, periodically
distributed piezoelectric inclusions. We consider the equations of linear
elasticity coupled with the electrostatic equation, with various kinds of
electric boundary conditions. We derive the corresponding effective models
when the thickness $a$ of the plate and the characteristic dimension $%
\varepsilon $ of the inclusions tend\ together to zero, in the two following
situations: first when $a \simeq \varepsilon$, second when $a/\varepsilon$
tends to zero.

\bigskip

{\bf \noindent Mathematics Subject Classifications (1991):}{\rm {\bf \ }}%
35B27 Homogenization of partial differential equations in media with
periodic structures, 73B27 nonhomogeneous materials and homogenization.

\bigskip

{\bf \noindent Key words:} linear elasticity, piezoelectricity,
homogenization, plate theory, composite materials, prescribed electric
potential, local electric circuits, nonlocal electric circuits, transfinite
networks, smart materials.

\newpage

\bigskip \noindent {\bf 1. Introduction }

\medskip \noindent {\bf 1.1 General. } This paper is the third and last part
of a systematic work devoted to the derivation of effective models for thin
piezoelectric/elastic composite plates\ including elementary electric
circuits connected to the upper and lower faces of piezoelectric
transducers. It is motivated by an important development of piezocomposites
used for instance for distributed control in vibroacoustics [8, 11 ,12, 25]
or as sensors in phased arrays. In [5], we considered three-dimensional
elastic plates with a small number of piezoelectric inclusions and derived
effective models when the thickness $a$ of the plate tends to zero. In [7],
effective models of thin plates with a large number of $\varepsilon $%
-periodically distributed piezoelectric inclusions have been obtained in the
case $\varepsilon <<a$ by letting $a$, $\varepsilon $ and $\varepsilon /a$
simultaneously tend to zero. The aim of the present paper is to consider the
two other possible asymptotic cases, namely $a/{\varepsilon \rightarrow 0}$
and $\varepsilon /a\rightarrow 1$. We note that the models for $%
a/\varepsilon \rightarrow 0$\ were already presented in the note [6].

\medskip \noindent As in [5, 7], different kinds of boundary conditions are
considered on the metallized upper faces of the inclusions, corresponding to
different possible types of control: prescribed electric potential (or
Dirichlet conditions) if the tension is controlled, prescribed electric
displacement field (or Neumann conditions) if the current is controlled,
local and nonlocal mixed conditions if the inclusions are connected to R-L-C
circuits (the nonlocal conditions corresponding to inclusions that are
connected to each other via R-L-C circuits).

\medskip \noindent Following a principle used in [5, 7], the derivation of
the models is made in the space of gradients of solutions. This\ leads to a\
more synthetic and readable presentation of the results. We combine the
two-scale convergence [1, 19] for homogenization and classical arguments of
thin plates theory [9, 10, 20]. Let us quote the pioneer work of Caillerie
[4], who considered the case of thin static elastic plates with periodic
coefficients, using the Tartar's method of oscillating test functions [2,
22]. However, in [4] the parameters $a$ and $\varepsilon $ tend successively
and independently to zero, except for $a\simeq \varepsilon $.

\medskip \noindent Despite their relative formal complexity, the effective
models have a rather simple structure. In the case $a/\varepsilon
\rightarrow 0$, for Dirichlet, and local mixed conditions (from a
mathematical point of view the prescribed electric displacement field
conditions happen to be a special case of mixed conditions and is not
treated separately), the limit model has the same form as the elastic plate
model, the influence of piezoelectric inclusion only appears in the
definition of the effective coefficients and as a source term. For nonlocal
mixed conditions, a coupling arises between mechanical effects and the
transverse component of the electric field, because of a Laplace operator in
the in-plane direction induced by the electric circuits. The resulting
operator is a special case of those encountered in the homogenization of
periodic electric circuits [15-17] or in transfinite networks [27, 28], the
later being analyzed from a different point of view. We also quote [21]
where a wide variety of in-plane operators are generated by a periodic
network of resistances. In the case $a/\varepsilon \rightarrow 1$, for
Dirichlet conditions the limit model also has the same form as the purely
elastic thin plate model, but a coupling arises even in the case of local
mixed conditions.

{\bf \medskip \noindent }Concerning thin piezoelectric structures, let us
also mention: Ghergu and al [14] who consider perforated piezoelectric
shells with fixed thickness; Licht and Weller [18, 26]; Sene [23] and
Figueiredo and Leal [13] who consider piezoelectric plates without
homogenization. Remark however that in [13] general {\it an}isotropic models
are considered. Finally, let us remark that the models used in
vibro-acoustic applications, as in [8, 11, 12, 25], are often based on Bloch
wave decompositions which seems a priori not compatible with the
homogenization method used in this paper. However, the homogenized model for
the wave equation in [3] builds a bridge between the two views and
constitutes a perspective for further works.

\medskip \noindent {\bf 1.2 Detailed contents. }Section 2 is devoted to the
setting of the initial 3-dimensional equations of static linearized
elasticity and piezoelectricity. The piezoelectric inclusions are assumed to
be strictly included in an insulating elastic matrix.

\medskip \noindent As our work is more about introducing piezoelectric
plates with elementary electric devices and about mixing homogenization and
plates theory, for simplicity and efficiency, following [4], we assume that
the material coefficients are constant in the thickness direction. This is
usually the case in applications, as the matrix and the piezoelectric
ceramics are homogeneous materials. In the same spirit, the upper and lower
faces of the piezoelectric inclusions are assumed to be metallized, that is
to be covered with a thin conductive film. However, from the mathematical
point of view, it might be interesting to obtain more general models by
removing these technical assumptions and considering fully non-homogeneous
materials as in [13], or multi-layered plate by adapting to the present work
to formalism proposed in [5].

\medskip \noindent The mechanical boundary conditions applied to the plate
are prescribed surface forces on its lower and upper faces and on part of
its lateral boundary, and prescribed mechanical displacement on the
remaining part. For the Maxwell-Gauss equations, we consider prescribed
electrical potential on the lower faces of the inclusions (in practice these
faces are connected to ground and the electric potential is zero). As
already mentioned in Section 1.1 various boundary conditions are considered
on the upper faces of the inclusions, corresponding to connections to
electric or electronic devices. These conditions are detailed in section
2.4. Some of them are, to our knowledge, unusual in plates theory, and thus,
constitute one of the interesting point of our work.

\medskip \noindent The weak formulations of the system are stated in Section
3. For a concise formulation covering all kinds of boundary conditions, we
adopt synthetic tensorial notations rather than fully extended formulae. We
strongly believe that this allows a better legibility of computations as
well as of limit models.

\medskip\noindent The precise assumptions on the data are presented in
Section 4. In particular, we give the correct scalings, or, from a more
concrete point of view, how electrical circuits have to be chosen to obtain
a significant influence on the effective behaviour of the material.
Resulting a priori estimates and first convergence results are given in
Sections 4.2 and 4.3.

\medskip \noindent Sections 5 and 6 are devoted to the statement of the main
results i.e. the effective two-dimensional plate model for each type of
electrical boundary condition in the case where the plate thickness is much
smaller than the inclusions size (Theorem 5.1, Section 5) and in the case
where these two parameters are of the same order (Theorem 6.1, Section 6).

\medskip \noindent Theorem 5.1 is proved in Section 7. The proof is in three
steps. The first one, which is mathematically the most difficult consists in
characterizing the two-scale limits of the strain and of the electrical
field. These results are new, even in the case of pure elasticity.In
Caillerie [4] - as two-scale convergence did not exist at the time - only
weak limits were considered. The second step consists in eliminating the
microscopic variable by computing the microscopic fields in terms of the
macroscopic fields. We use here the classical arguments of linear
homogenization. The third step consists in eliminating the transverse
components, or part of the transverse components of the fields (according to
the model), that are computed with respect to the other ones. This
elimination slightly departs from the classical plates theory, because of
the non standard boundary conditions on the upper and lower faces of the
inclusions that are considered in the present paper.

\smallskip \noindent Theorem 6.1 is proved in Section 8. As here $a\simeq
\varepsilon $ is assumed, the proof is only in two steps. First,
characterization of the limit, second, simultaneous elimination of the local
and of (part of) the transverse components.

\smallskip \noindent We use the same formalism as in [5,7], based on
tensorial notations and products, and on simple algebraic operations such as
projections. It allows to deal relatively easily with complex computations.
Completely explicit formulae would require a lot of room, to the detriment
of legibility. Step 2 and 3 of our proofs are almost formal computation and
are easily adapted from one variant to the other. A coupling with the
formalism introduced in [5] for multilayered plates models is easily
conceivable.

\bigskip \noindent {\bf 2. Equations of 3-dimensional piezoelectricity}

\medskip \noindent {\bf 2.1 Geometry}. The three-dimensional plate with
thickness $a>0$ is represented by ${\Omega }^{a}={\ \omega }\times ]-a,a[$, $%
{\omega }$ being a bounded domain of ${\mathbb R}^{2}$, see Figure
1. Using the change of scales and variables introduced in [10], we
shall work on the fixed domain ${\Omega }={\omega }\times ]-1,1[$.
The domain ${\omega }$ is
divided into two subdomains ${\omega }_{1}^{{\varepsilon }}$ and ${\omega }%
_{2}^{{\varepsilon }}$ that are constructed as follows.\ Let $%
Y=]-1/2,1/2[^{2}$ and $Y_{1}$ be a strict smooth subdomain of $Y$, see
Figure 2, let ${\mathbb I}^{\varepsilon }$ be the set of multi-index ${\bf i}%
=(i_{1},i_{2})\in {\mathbb Z}^{2}$ such that ${\varepsilon (}{\bf
i+}Y_{1})$ is strictly included in ${\omega }$; then $\omega
_{1}^{\varepsilon
}=\bigcup\limits_{{\bf i}\in {\mathbb I}^{\varepsilon }}{\varepsilon (}{\bf i+}%
Y_{1})$ while ${\omega }_{2}^{{\varepsilon }}={\omega }\setminus \overline{%
\omega }_{1}^{\varepsilon }$. Let $b=(a,{\varepsilon })$; the set ${\Omega }%
_{2}^{b}={\omega }_{2}^{{\varepsilon }}\times ]-a,a[$ represents the elastic
matrix of the body, while ${\Omega }_{1}^{b}={\omega }_{1}^{{\varepsilon }%
}\times ]-a,a[$ is the set of all piezoelectric inclusions.

\begin{figure}[htbp]
\begin{center}
\includegraphics[width=16 cm]{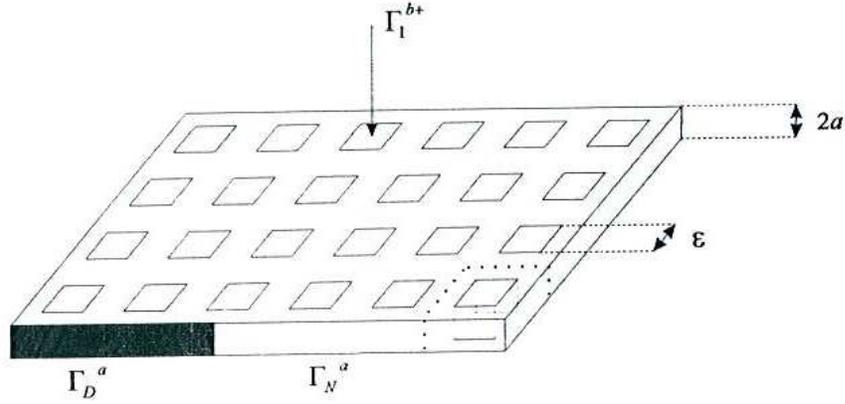}
\end{center}
\caption{\text{ Composite plate with piezoelectric inclusions.} }
\end{figure}

\noindent The boundary of ${\omega }$ is divided into two regular parts $%
\gamma _{D}$ and $\gamma _{N}$, with $|\gamma _{D}|>0$. The boundary of ${%
\Omega }^{a}$ is thereby divided into ${\Gamma }_{D}^{a}={\gamma }_{D}\times
]-a,a[$ and ${\Gamma }_{N}^{a}$ = $({\gamma }_{N}\times ]-a,a[)$ $\cup $ $({%
\omega }\times \{-a,a\})$. The boundary of ${\Omega }_{1}^{b}$ is divided
into ${\Gamma }_{1}^{b+}={\omega }_{1}^{{\varepsilon }}\times \{a\}$, ${%
\Gamma }_{1}^{b-}={\omega }_{1}^{{\varepsilon }}\times \{-a\}$ and ${\Gamma }%
_{1}^{b}$ = ${\partial }{\omega }_{1}^{\varepsilon }\times ]-a,a[$. The
outer unit normals to the boundaries of $\Omega ^{a}$ and $Y$ are denoted by
${\bf n}$ and ${\bf n}_{Y}$, respectively.

\medskip \noindent For any inclusion ${\varepsilon (}{\bf i+}Y_{1})\times
]-a,a[$ such that ${\bf i\in }{\mathbb I}^{\varepsilon }$, the mean
value on
the upper face ${\varepsilon (}{\bf i+}Y_{1})\times \{a\}$ is denoted by $%
<.>_{{\bf i}}.$ For any function $\psi $ on $\Omega ^{a}$, ${\psi }_{{\bf i}%
} $ designates its restriction to the inclusion ${\varepsilon (}{\bf i+}%
Y_{1})\times ]-a,a[$.

\smallskip \noindent The space variables are $x^{a}=(\hat{x},x_{3}^{a})\in
\Omega ^{a}$ where $x_{3}^{a}\in ]-a,a[$, $\hat{x}=(x_{1},x_{2})\in {\omega }
$ and ${y}=(y_{1},y_{2})\in Y$. The derivatives with respect to $x_{\alpha }$%
, $x_{3}^{a}$ and $y_{\alpha }$ are denoted by ${\partial }_{\alpha }$, ${%
\partial }_{3}$ and ${\partial }_{y_{\alpha }}$, respectively.

\smallskip \noindent When referring to the fixed domain $\Omega $, the
geometric notation is the same, the subscript $a$ being removed if necessary.

\begin{figure}[htbp]
\begin{center}
\includegraphics[keepaspectratio=true, height=5.5cm]{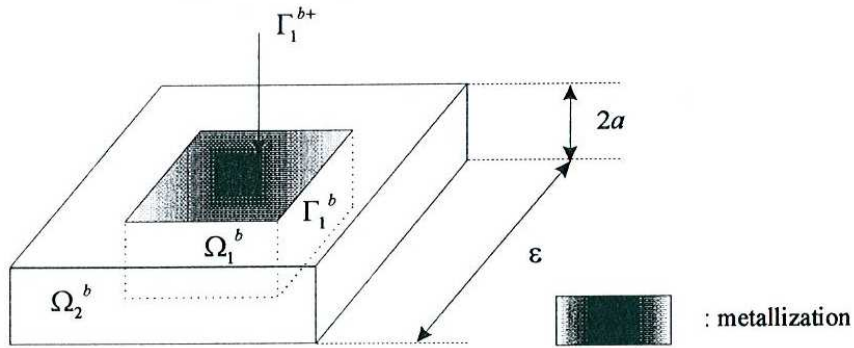}
\end{center}
\caption{\text{ Elementary cell, with piezoelectric inclusion and
metallization.} }
\end{figure}

\medskip \noindent {\bf 2.2 Other notations. } Bold characters are used for
vector and matrix valued functions and for the corresponding functional
spaces. We constantly use Einstein's convention of summation on repeated
indices with the additional convention that latin and greek indices are
varying from 1 to 3 and from 1 to 2 respectively. Throughout the paper $c$
and $C$ designate generic positive constants, not depending on $a$ and $%
\varepsilon $.

\medskip \noindent Also:
\begin{equation}
inutile
\label{1003}
\end{equation}
Last, vector and vectorial notations are a priori meant $in$ $line$.
When a
given {\em vector }${\bf V}$ is meant {\em in column}, we write:{\em \ }$^{t}%
{\bf V}$. For instance for
\[
{\bf V}=\left( \left( s_{{\alpha }\beta }({\bf v})\right) _{\alpha ,\beta
=1,2},\text{ }\left( a^{-1}\,s_{\alpha 3}({\bf v})\right) _{\alpha =1,2},\
a^{-2}\,s_{33}({\bf v})\right) \text{,}
\]%
with $^{t}{\bf V}$, we mean%
\[
^{t}{\bf V=}\left( {\bf
\begin{tabular}{l}
$\left( s_{{\alpha }\beta }(v)\right) _{\alpha ,\beta =1,2}$ \\
$\left( a^{-2}\,s_{\alpha 3}(v)\right) _{\alpha =1,2}$ \\
$a^{-2}\,s_{33}(v)$%
\end{tabular}%
}\right) .
\]%
Note that this convention leads us to writings that slightly differ from the
one in [6, 7], but it seems more coherent to us now.

\medskip \noindent {\bf 2.3 Equations of 3-dimensional piezoelectricity.}
The mechanical displacements ${\bf u}^{b}=(u_{i}^{b})_{i=1,2,3}$ and the
electric potential $\varphi ^{b}$ are governed by the linearized equations
of piezoelectricity in their static version,

\begin{equation}
\left\{
\begin{array}{l}
-{\partial }_{j}{\sigma }_{ij}^{b}=f_{i}^{b}\text{ in }{\Omega }^{a}\text{, }%
\,\,\,{\sigma }_{ij}^{b}n_{j}=g_{i}^{b}\text{ on }{\Gamma }%
_{N}^{a},\,\,\,u_{i}^{b}=0\text{ on }{\Gamma }_{D}^{a},\medskip \\
-{\partial }_{j}D_{j}^{b}=0\text{ in }{\Omega }_{1}^{b},\,\,D_{j}^{b}.n_{i}=0%
\text{ on }{\Gamma }_{1}^{b},%
\end{array}%
\right.  \label{1.b}
\end{equation}%
for i = 1, 2, 3, where

\begin{equation}
\left\{
\begin{array}{l}
{\sigma }_{ij}^{b}=R_{ijkl}^{\varepsilon }s_{kl}({\bf u}^{b})+d_{kij}^{%
\varepsilon }{\partial }_{k}\varphi ^{b},\medskip \\
D_{j}^{b}=-d_{jkl}^{\varepsilon }s_{kl}({\bf u}^{b})+c_{jk}^{\varepsilon }{%
\partial }_{k}\varphi ^{b},%
\end{array}%
\right.  \label{60}
\end{equation}%
\begin{equation}
\forall \,{\bf v}\in {\bf H}^{1}({\Omega }^{a}),\;\;\;s_{ij}({\bf v})={\frac{%
1}{2}}({\partial }_{i}v_{j}+{\partial }_{j}v_{i}).  \label{1005}
\end{equation}%
\noindent The assumptions on the volume and surface forces ${\bf f}%
^{b}:=(f_{i}^{b})_{i=1,2,3}$ and ${\bf g}^{b}:=(g_{i}^{b})_{i=1,2,3}$ are
specified below, in Section 4.1. The stiffness tensor ${\bf R}^{\varepsilon
}:=(R_{ijkl}^{{\varepsilon }})_{i,j,k,l=1,2,3}$, the piezoelectricity tensor
${\bf d}^{\varepsilon }:=(d_{kij}^{\varepsilon })_{i,j,k=1,2,3}$ and the
permitivity tensor ${\bf c}^{\varepsilon }:=(c_{ij}^{\varepsilon
})_{i,j=1,2,3}$ are assume to satisfy the symmetry conditions:%
\begin{equation}
\forall i,j,k,l\in \{1,2,3\},\;\;\;R_{ijkl}^{{\varepsilon }}=R_{klij}^{{%
\varepsilon }}=R_{jikl}^{{\varepsilon }},\;\;c_{ij}^{\varepsilon
}=c_{ji}^{\varepsilon },\;\;d_{kij}^{\varepsilon }=d_{kji}^{\varepsilon }%
\text{.}  \label{2}
\end{equation}%
\noindent Note that, as we assume that the inclusions are electrically
insulated from the elastic matrix, no Gauss-Maxwell equation is needed in ${%
\Omega}_{2}^{b}$. However, we let for convenience
\begin{equation}
{\bf c}^{\varepsilon }={\bf 0},\;\;{\bf d}^{\varepsilon }={\bf 0}\;\;\;\text{%
in }\Omega _{2}^{b}\text{.}  \label{1001}
\end{equation}%
We go now into detail about the different kinds of electric boundary
conditions that are considered in this paper.

\medskip \noindent {\bf 2.4 Electric boundary conditions on }$\Gamma
_{1}^{b+}{\bf \cup }\Gamma _{1}^{b-}${\bf .} The four kinds of electric
boundary conditions are summarized here. Further explanations and comments
can be found in [7].

\medskip \noindent (i) {\em Prescribed electric potential (Dirichlet
conditions):}%
\begin{equation}
\varphi ^{b}=\varphi _{m}^{b}+a\varphi _{c}^{b}\text{ on }{\Gamma }%
_{1}^{b+},\,\ \ \ \,\varphi ^{b}=\varphi _{m}^{b}-a\varphi _{c}^{b}\text{ on
}{\Gamma }_{1}^{b-},  \label{39}
\end{equation}%
$\varphi _{m}^{b}$ and $\varphi _{c}^{b}$ being two given functions on $%
\omega _{1}^{\varepsilon }$.

\medskip \noindent (ii) {\em Prescribed electric displacement field (Neumann
conditions), local and nonlocal electric circuits (local and nonlocal mixed
conditions): }

\medskip \noindent Let us introduce the shift operators defined from ${\mathbb I%
}^{\varepsilon }$ to ${\mathbb Z}^{2}$ by%
\[
\left\{
\begin{array}{l}
T_{+1}^{1}:{\bf i\mapsto }(i_{1}+1,i_{2}),\text{\ }T_{+1}^{2}:{\bf i\mapsto }%
(i_{1},i_{2}+1), \\
T_{-1}^{1}:{\bf i\mapsto }(i_{1}-1,i_{2}),\;T_{-1}^{2}:{\bf i\mapsto }%
(i_{1},i_{2}-1).%
\end{array}%
\right.
\]%
Let $\varphi _{m}^{b}$ be a given function on $\omega _{1}^{\varepsilon }$.
For all ${\bf i}\in {\mathbb I}^{\varepsilon }$, we pose%
\[
\bar{\varphi}_{{\bf i}}^{b}=\varphi _{{\bf i}}^{b}-\varphi _{m,{\bf i}}^{b}.
\]%
The boundary conditions for the electrostatic equation (\ref{1.b})$_{2}$ on $%
{\Gamma }_{1}^{b+}\cup {\Gamma }_{1}^{b-}$ are then%
\begin{equation}
\left\{
\begin{array}{l}
\forall {\bf i}\in {\mathbb I}^{\varepsilon },\;\;<{\bf D}^{b}.{\bf n}>_{{\bf i}%
}=\dfrac{G_{1}}{{a{\varepsilon }^{2}}}\sum\limits_{\alpha =1}^{2}(\bar{%
\varphi}_{T_{-1}^{\alpha }({\bf i})}^{b}-2\bar{\varphi}_{{\bf i}}^{b}+\bar{%
\varphi}_{T_{+1}^{\alpha }({\bf i})}^{b})-{\dfrac{G}{a}}\bar{\varphi}_{{\bf i%
}}^{b}+h^{b}\text{\ \ \ on}\,{\Gamma }_{1}^{b+}, \\
\varphi ^{b}=\varphi _{m}^{b}\text{\ \ \ on }{\Gamma }_{1}^{b-}.%
\end{array}%
\right.  \label{41}
\end{equation}%
The finite difference operator on $\bar{\varphi}_{{\bf i}}^{b}$ is completed
by the analog of a discrete Neumann boundary condition to the free ends of
the circuit
\begin{equation}
\bar{\varphi}_{T_{-1}^{\alpha }({\bf i})}^{b}-\bar{\varphi}_{{\bf i}}^{b}=0%
\text{\ \ if\ \ }T_{-1}^{\alpha }({\bf i})\notin {\mathbb I}^{\varepsilon }%
\text{,\ \ \ \ }\bar{\varphi}_{T_{+1}^{\alpha }({\bf i})}^{b}-\bar{\varphi}_{%
{\bf i}}^{b}=0\text{\ \ if\ \ }T_{+1}^{\alpha }({\bf i})\notin {\mathbb I}%
^{\varepsilon }  \label{41-1}
\end{equation}%
which also means that the current is vanishing in the corresponding branches.

\medskip \noindent In (\ref{41}), $G_{1}$ and $G$ are given nonnegative
constants. If $G$ $=$ $G_{1}=0$ the conditions on ${\Gamma }_{1}^{b+}$ are
Neumann condition (prescribed electric displacement field). If $G$ $>0$ and $%
G_{1}=0$ those are local mixed conditions. If $G$ $G_{1}>0$ those are
nonlocal mixed conditions. If $G_{1}=0 $, the condition (\ref{41-1}) is of
course not relevant.

\medskip \noindent Although from a physical point of view it has its own
meaning, the case of Neumann conditions does not need to be treated
separately; one just have to let $G=0$ in the effective models obtained in
Theorem 5.1 and Theorem 6.1 for local mixed conditions. In the situation $%
a/\varepsilon \rightarrow 0$, local and nonlocal conditions lead to
different developments and limit models.

\medskip \noindent As we assume all the faces to be metallized, the
functions $\bar{\varphi}^{b}$, $\varphi _{m}^{b}$\ and $\bar{\varphi}%
_{c}^{b} $\ are constant on each metallized face of inclusions and assuming
that the current is provided by a single wire, the same holds true for $%
h^{b} $.

\medskip \noindent Finally, in order to use, as much as possible, common
formulations for the different boundary conditions, we define $h^{b}$, $%
\varphi _{c}^{b}$, $G$\ and $G_{1}$\ in all cases, with the conventions that
\begin{equation}
h^{b}=0\;\text{and }G=G_{1}=0\text{ for Dirichlet conditions,}  \label{492}
\end{equation}%
and%
\begin{equation}
\varphi _{c}^{b}=0\;\text{for mixed conditions.}  \label{493}
\end{equation}

\bigskip \noindent {\bf 3. Weak formulations }

\medskip \noindent The aim of the present section is the statement of the
weak formulation of the equations of Section 2, on the fixed domain ${\Omega
}$. We use the standard change of variables $x^{a}\rightarrow
x=(x_{1},x_{2},x_{3}^{a}{/}a)$ and the appropriate scaling for volume
forces, surface forces and displacements fields, see [10] and also [9],%
\[
\left\{
\begin{tabular}{l}
${\bf \forall }x\in \Omega $,$\ \ {\bf \hat{u}}%
^{b}(x)=(u_{1}^{b}(x^{a}),u_{2}^{b}(x^{a}),au_{3}^{b}(x^{a})),$ \\
${\bf \forall }x\in \Omega $,$\ \ {\bf \hat{f}}%
^{b}(x)=(f_{1}^{b}(x^{a}),f_{2}^{b}(x^{a}),a^{-1}f_{3}^{b}(x^{a})),$ \\
${\bf \forall }x\in \gamma _{N}\times ]-1,1[$,$\ \ {\bf \hat{g}}%
^{b}(x)=(g_{1}^{b}(x^{a}),g_{2}^{b}(x^{a}),a^{-1}g_{3}^{b}(x^{a})),$ \\
${\bf \forall }x\in \omega \times \left\{ -1,1\right\} $,$\ \ {\bf \hat{g}}%
^{b}(x)=a^{-1}(g_{1}^{b}(x^{a}),g_{2}^{b}(x^{a}),a^{-1}g_{3}^{b}(x^{a})).$%
\end{tabular}%
\right.
\]%
The current source $h^{b}$, the electric potential $\varphi ^{b}$, $\varphi
_{m}^{b}$ and $\varphi _{c}^{b}$ are left unchanged. As in the sequel, we
only work on the reference domain $\Omega $, no confusion might occur, so
for simplicity we keep the notation ${\bf u}^{b}$, ${\bf f}^{b}$, ${\bf g}%
^{b}$, $h^{b},$ $\varphi ^{b}$, $\varphi _{m}^{b}$ and $\varphi _{c}^{b}$,
without hats.

\smallskip \noindent For any ${\bf V=}({\bf v},\psi {\bf )}\in {\bf H}^{1}({%
\Omega })\times H^{1}(\Omega _{1}^{\varepsilon }),$ we let
\begin{equation}
\left\{
\begin{tabular}{l}
${\bf K}^{a}({\bf v})=\left( K_{{\alpha }\beta }^{a}({\bf v}),K_{{\alpha }%
3}^{a}({\bf v})),{K}_{33}^{a}({\bf v})\right) =\left( s_{{\alpha }\beta }(%
{\bf v}),a^{-1}\,s_{\alpha 3}({\bf v}),a^{-2}\,s_{33}({\bf v})\right)
,\medskip $ \\
${\bf L}^{a}(\psi )=\left( L_{\alpha }^{a}(\psi ),L_{3}^{a}(\psi )\right)
=\left( {\partial }_{\alpha }\psi ,a^{-1}\,{\partial }_{3}\psi \right)
,\medskip $ \\
${\bf M}^{a}{\bf (V})=\left( {\bf K}^{a}({\bf v}),{\bf L}^{a}(\psi )\right)
. $%
\end{tabular}%
\right.  \label{42}
\end{equation}%
\noindent We put together the tensors ${\bf R}^{\varepsilon }$, ${\bf d}%
^{\varepsilon }$ and ${\bf c}^{\varepsilon }$ in a global
stiffness-piezoelectricity-permitivity tensor ${\cal R}^{\varepsilon }$,
which is the $10\times 10$ symmetric matrix written in a format compatible
with (\ref{42}):

\begin{equation}
{\cal R}^{\varepsilon }=%
\begin{array}{r}
\left(
\begin{array}{ccccc}
R_{{\alpha }\beta {\gamma }{\delta }}^{\varepsilon } & 2R_{{\alpha }\beta {%
\gamma }3}^{\varepsilon } & R_{{\alpha }\beta 33}^{\varepsilon } & d_{{%
\gamma }{\alpha }\beta }^{\varepsilon } & d_{3{\alpha }\beta }^{\varepsilon }
\\
2R_{{\alpha }3{\gamma }{\delta }}^{\varepsilon } & 4R_{{\alpha }3{\gamma }%
3}^{\varepsilon } & 2R_{{\alpha }333}^{\varepsilon } & 2d_{{\gamma }{\alpha }%
3}^{\varepsilon } & 2d_{3{\alpha }3}^{\varepsilon } \\
R_{33{\gamma }{\delta }}^{\varepsilon } & 2R_{33{\gamma }3}^{\varepsilon } &
R_{3333}^{\varepsilon } & d_{{\gamma }33}^{\varepsilon } &
d_{333}^{\varepsilon } \\
-d_{{\alpha }{\gamma }{\delta }}^{\varepsilon } & -2d_{{\alpha }{\gamma }%
3}^{\varepsilon } & -d_{{\alpha }33}^{\varepsilon } & c_{{\alpha }{\gamma }%
}^{\varepsilon } & c_{{\alpha }3}^{\varepsilon } \\
-d_{3{\gamma }{\delta }}^{\varepsilon } & -2d_{3{\gamma }3}^{\varepsilon } &
-d_{333}^{\varepsilon } & c_{3{\gamma }}^{\varepsilon } &
c_{33}^{\varepsilon }%
\end{array}%
\right) .%
\end{array}
\label{tenseur}
\end{equation}%
The linear forms associated with the mechanical load and the electrical
current source are
\[
\ell _{u}^{b}({\bf v})=\int\nolimits_{{\Omega }}f_{i}^{b}\,v_{i}\,\text{d}%
x+\int\nolimits_{{\Gamma }_{N}}g_{i}^{b}\,v_{i}\,\text{d}s,\text{ and }%
\,\,\ell _{\varphi }^{b}(\tilde{L}_{3})=\int\nolimits_{{\Omega }%
_{1}^{\varepsilon }}h^{b}\tilde{L}_{3}\,\text{d}x.
\]%
\noindent Note that $h^{b}$ is a priori defined on $\Gamma _{1}^{\varepsilon
+}$ (or equivalently on $\omega_{1}^{\varepsilon}$) but is trivially
extended to a function on $\Omega_{1}^{\varepsilon}$ which does not depend
on $x_3$. \noindent Given the assumption on metallization, the set of
admissible electric potential is
\[
H_{c}^{1}({\Omega }_{1}^{\varepsilon })=\{\psi \in H^{1}({\Omega }%
_{1}^{\varepsilon });\psi \text{ is constant on each connected part of }%
\Gamma _{1}^{\varepsilon +}\cup \Gamma _{1}^{\varepsilon -}\}.
\]%
\noindent The relevant functional space is then specific to each electric
boundary condition.

\smallskip {\em \noindent (i) Dirichlet conditions:}%
\[
\left\{
\begin{tabular}{l}
$%
\begin{array}{l}
{\bf W}_{D}^{b}={\bf W}^{\varepsilon }(\varphi _{m}^{b},\varphi _{c}^{b})=\{(%
{\bf v},\,\varphi )\in {\bf H}^{1}(\Omega )\times H_{c}^{1}(\Omega
_{1}^{\varepsilon });\,{\bf v}={\bf 0}\text{ on }\Gamma _{D}, \\
\text{{\bf \ \ \ \ \ \ \ \ \ \ \ \ \ \ \ \ \ \ \ }}\varphi =\varphi
_{m}^{b}+a\varphi _{c}^{b}\text{ on }\Gamma _{1}^{\varepsilon +},\;\varphi
=\varphi _{m}^{b}-a\varphi _{c}^{b}\text{ on }\Gamma _{1}^{\varepsilon -}\},%
\end{array}%
$ \\
${\bf W}^{\varepsilon }={\bf W}^{\varepsilon }(0,0).$%
\end{tabular}%
\right.
\]%
{\em (ii) Mixed conditions:}%
\[
\left\{
\begin{array}{l}
{\bf W}_{D}^{b}={\bf W}^{\varepsilon }(\varphi _{m}^{b})=\{({\bf v}%
,\,\varphi )\in {\bf H}^{1}(\Omega )\times H_{c}^{1}(\Omega
_{1}^{\varepsilon });\,{\bf v}={\bf 0}\text{ on }\Gamma _{D},\,\varphi
=\varphi _{m}^{b}\text{ on }\Gamma _{1}^{\varepsilon -}\}, \\
{\bf W}^{\varepsilon }={\bf W}^{\varepsilon }(0).%
\end{array}%
\right.
\]%
\noindent The backward difference operator $\nabla _{\hat{x}}^{\varepsilon }$
is defined inclusion by inclusion by
\[
\forall {\bf i}\in {\mathbb I}^{\varepsilon },\;\forall \psi \in H_{c}^{1}({%
\Omega }_{1}^{\varepsilon }),\;\;(\nabla _{\hat{x}}^{\varepsilon }\psi )_{%
{\bf i}}=\varepsilon ^{-1}\text{ }(\psi _{{\bf i}}-\psi _{T_{-1}^{1}({\bf i}%
)},\,\psi _{{\bf i}}-\psi _{T_{-1}^{2}({\bf i})}).
\]%
\smallskip \noindent Letting ${\bf U}^{b}=({\bf u}^{b},\varphi ^{b})$, with
the conventions (\ref{1001}), (\ref{492}) and (\ref{493}), the weak
formulations on the scaled domain ${\Omega }$ for the coupled problems (\ref%
{1.b})-(\ref{39})-(\ref{41}), are summarized by:

\begin{equation}
\left\{
\begin{array}{l}
{\bf U}^{b}=({\bf u}^{b},\varphi ^{b})\in {\bf W}_{D}^{b}\text{,}\;\ \ \text{%
and for all \ }{\bf V}=({\bf v},{\psi })\in {\bf W}^{\varepsilon }, \\
\ \ \ \ \dint\nolimits_{\Omega }\,{\bf M}^{a}({\bf V})\ {\cal R}%
^{\varepsilon }\ ^{t}{\bf M}^{a}({\bf U}^{b})\,\text{d}x+{2}\dint\nolimits_{{%
\Omega }_{1}^{\varepsilon }}G{\cal M}(L_{3}^{a}(\varphi ^{b})){\cal M}%
(L_{3}^{a}({\psi }))\,\text{d}x \\
\;\;\;\;\;\;\;\;\;\;\;\;\;+2\dint\nolimits_{{\Omega }_{1}^{\varepsilon
}}G_{1}\nabla _{\hat{x}}^{\varepsilon }{\cal M}(L_{3}^{a}(\varphi ^{b})){\bf %
.}\nabla _{\hat{x}}^{\varepsilon }{\cal M}(L_{3}^{a}({\psi }))\,\text{d}%
x=\ell _{u}^{b}({\bf v})+\ell _{\varphi }^{b}(L_{3}^{a}({\psi })).%
\end{array}%
\right.  \label{63}
\end{equation}

\bigskip \noindent {\bf 4. Assumptions on the data - A priori estimates -
Convergences }

\medskip \noindent {\bf 4.1 Two-scale convergence.} As two-scale convergence
is an important tool of the paper, before stating the assumptions on the
data and the first convergence results, that are expressed in terms of
two-scale convergence, let us recall some basic facts about it. Let $%
C_{\sharp }^{\infty }(Y)$ designate the space of $C^{\infty } $
functions on ${\mathbb R}^{n}$ that are $Y-$periodic.

\medskip \noindent {\bf D\'{e}finition 1 }(Allaire [1]) Let $(u^{\varepsilon
})_{\varepsilon >0}$ be a family of $L^{2}(\omega )$ and $u\in L^{2}(\omega
\times Y)$. We say that $(u^{\varepsilon })_{\varepsilon >0}$ two-scale
converges to $u$ if for any $v\in {\cal D}(\omega ;C_{\sharp }^{\infty }(Y))$%
\[
\lim\limits_{\varepsilon \rightarrow 0}\int_{\omega }u^{\varepsilon }(\hat{x}%
)v(\hat{x},\hat{x}/\varepsilon )\text{d}\hat{x}=\int_{\omega \times Y}u(\hat{%
x},y)v(\hat{x},y)\text{d}\hat{x}\text{d}y.
\]

\medskip \noindent The important fact here is that for any bounded family in
$L^{2}(\omega )$, there is a subsequence that two-scale converges to some
limit $u$. Two-scale convergence is a more accurate notion than usual weak
convergence, in the sense that the two-scale converging family also weakly
converges to $\dint_{Y}u(.,y)$d$y$.

\medskip \noindent Since we need two-scale convergence for functions defined
on $\omega _{1}^{\varepsilon }$, we also use the following practical
definition.

\medskip \noindent {\bf D\'{e}finition 2 \ }A family $(u^{\varepsilon
})_{\varepsilon >0}$ of $L^{2}(\omega _{1}^{\varepsilon })$ is said to
two-scale converge to a limit $u$ in $L^{2}(\omega \times Y_{1})$ if $u\in
L^{2}(\omega \times Y_{1})$ and $\left( P^{\varepsilon }u^{\varepsilon
}\right) $ two-scale converges to $Pu$ in the sense of definition 1, where $%
P^{\varepsilon }$ and $P$ designate the extension by $0$ of functions on $%
\omega _{1}^{\varepsilon }$ to functions on $\omega $ and of functions on $%
\omega \times Y_{1}$ to functions $\omega \times Y$ respectively.

\medskip \noindent In our problems, there are no oscillations in the $x_{3}$%
-direction. Still, the above definitions evidently apply with $x_{3}$ as a
dummy variable (as it is the case for the time variables $t$ in other
contexts). So in what follows when writing that $(u^{\varepsilon })$ of $%
L^{2}(\Omega )$ two-scale converges to $u$ in $L^{2}(\Omega \times Y)$, we
mean that for any $v\in {\cal D}(\Omega ;C_{\sharp }^{\infty }(Y))$,%
\[
\lim\limits_{\varepsilon \rightarrow 0}\int_{\Omega }u^{\varepsilon }(x)v(x,%
\hat{x}/\varepsilon )\text{d}x=\int_{\Omega \times Y}u(x,y)v(x,y)\text{d}x%
\text{d}y,
\]

\medskip \noindent and similarly for $(u^{\varepsilon})$ in $L^{2}(\Omega
\times Y_{1})$.

\medskip \noindent Remark that the convergences in Definitions 1 and 2 are
weak convergences. In general, the two-scale limit of a product is not the
product of the two-scale limits. An additional assumption is needed. Let $%
(u^{\varepsilon })$ and $(v^{\varepsilon })$ be two two-scale converging
family in $L^{2}(\Omega )$ with limits $u$ and $v$ respectively. Then, if in
addition
\begin{equation}
\lim\limits_{\varepsilon \rightarrow 0}\left\Vert u^{\varepsilon
}\right\Vert _{L^{2}(\omega )}=\left\Vert u\right\Vert _{L^{2}(\omega \times
Y)},  \label{strong two scale}
\end{equation}%
the following statement hold true for any regular test functions $\varphi :$
\begin{equation}
\lim\limits_{\varepsilon \rightarrow 0}\int_{\omega }u^{\varepsilon }(\hat{x}%
)v^{\varepsilon }(\hat{x})\varphi (\hat{x},\hat{x}/\varepsilon )\text{d}\hat{%
x}=\int_{\omega \times Y}u(\hat{x},y)\varphi (\hat{x},y){\varphi }\text{d}%
\hat{x}\text{d}y.  \label{cv produit}
\end{equation}%
Assumption (\ref{strong two scale}) is a kind of strong two-scale
convergence notion. See for example Theorem 1.8 and the proof of Theorem 2.3
in Allaire [1] for a proof of (\ref{cv produit}).

\medskip \noindent {\bf 4.1 Assumptions on the data and first convergence
results. }The tensors ${\bf R}^{\varepsilon }$, ${\bf d}^{\varepsilon }$ and
${\bf c}^{\varepsilon }$ constituting the
stiffness-piezoelectricity-permitivity tensor ${\cal R}^{\varepsilon }$ are
assumed to satisfy (\ref{2}) and
\begin{equation}
\left\{
\begin{array}{l}
({\cal R}^{\varepsilon })\in {\bf L}^{\infty }(\Omega )\,\text{and two-scale
converges }\text{to some }{\cal R}\in {\bf L}^{\infty }({\Omega }\times
Y),\medskip \\
||{\cal R}^{\varepsilon }||_{{\bf L}^{\infty }({\Omega })}\leq C,\,\,{\cal R}%
^{\varepsilon }\,\,\text{does not depend on}\,\text{ }x_{3},\medskip \\
\forall \,{\bf K=(}K_{ij})\in {\mathbb R}^{3\times 3}\text{ with }%
K_{ij}=K_{ji},\;\;{\bf K}\ {\bf R}^{\varepsilon }\ ^{t}{\bf K}\geq c||{\bf K}%
||^{2}\;\;\text{a.e. in }\omega ,\medskip \\
\forall \,{\bf L}\in {\mathbb R}^{3},\;\;^{t}{\bf L}\ {\bf
c}^{\varepsilon }\
{\bf L}\geq c||{\bf L}||^{2}\;\;\text{a.e. in }\omega _{1}^{\varepsilon },%
\end{array}%
\right.  \label{8}
\end{equation}

\medskip \noindent and

\begin{equation}
{\lim\limits_{\varepsilon \rightarrow 0}}||{\cal R}^{\varepsilon }||_{{\bf L}%
^{2}({\omega })}=||{\cal R}||_{{\bf L}^{2}({\omega }\times Y)}.  \label{8bis}
\end{equation}

\medskip \noindent {\bf Remark:} As pointed above, assumption (\ref{8bis})
is a rather vague, but general assumption that allows to go to the limit in
product of two-scale converging functions; namely $\left( {\bf M}^{a}{\bf (U}%
^{b})\right) $ and $\left( {\cal R}^{\varepsilon }\right) $. It is of course
fulfilled by data satisfying classical periodicity conditions, for instance $%
{\cal R}^{\varepsilon }(\hat{x})={\cal R}(\hat{x},\hat{x}/\varepsilon )$
where ${\cal R}$ is a given tensor in ${\bf L}^{\infty }(\omega ;C_{\sharp
}^{{}}(Y))$, or also, which is more relevant in the present work, in ${\bf C}%
^{\infty }(\omega ;L_{\sharp }^{\infty }(Y))$. See Allaire [1] if necessary.

\medskip \noindent Coercivity for ${\bf c}^{\varepsilon }$ and ${\bf R}%
^{\varepsilon }$, together with the symmetry assumptions (\ref{2}) for ${\bf %
d}$, implies coercivity for ${\cal R}^{\varepsilon }$. Conversely, two-scale
convergence for ${\cal R}^{\varepsilon }$ implies two-scale convergence for $%
{\bf R}^{\varepsilon }$, ${\bf c}^{\varepsilon }$ and ${\bf d}^{\varepsilon
} $. The corresponding limits are denoted by ${\bf R}$, ${\bf c}$ and ${\bf d%
}$. The other data are assumed to satisfy
\begin{equation}
\left\{
\begin{array}{l}
{\bf f}^{b}\in {\bf L}^{2}({\Omega }),\text{ }{\bf g}^{b}\in {\bf H}^{1/2}({%
\Gamma }_{N}),\medskip \  \\
({\bf f}^{b})\text{ converges weakly in }{\bf L}^{2}(\Omega )\text{ to some
limit }{\bf f},\medskip \  \\
({\bf g}^{b})\text{ converges weakly in }{\bf L}^{2}({\Gamma }_{N})\,\,\text{%
to some limit }{\bf g}\text{,}%
\end{array}%
\right.  \label{9}
\end{equation}%
\begin{equation}
\left\{
\begin{array}{l}
h^{b},\,\,\varphi _{m}^{b}\,\,\text{and}\,\,\varphi _{c}^{b}\,\,\text{are
constant on each inclusion},\medskip \\
(h^{b})\text{ two scale-converges in }\,L^{2}({\omega }\times Y_{1})\,\,%
\text{ to some limit }h\in L^{2}({\omega }),\medskip \\
(\varphi _{m}^{b})\text{ two-scale converges in }\,L^{2}({\omega }\times
Y_{1})\,\,\text{ to some limit }\varphi _{m}\in H^{1}({\omega }),\medskip \\
(\varphi _{c}^{b})\text{ two-scale converges in }\,L^{2}({\omega }\times
Y_{1})\,\,\text{ to some limit }\varphi _{c}\in L^{2}({\omega }).%
\end{array}%
\right.  \label{49}
\end{equation}

\noindent We observe that, because $\varphi _{c}^{b}$, $\varphi _{m}^{b}$
and $h^{b}$ are constant on each inclusion, their two-scale limits do not
depend on $y$ in $Y_{1}$.

\medskip \noindent Last, let us introduce the space of Kirchhoff-Love's
displacement fields:
\begin{equation}
{\bf V}_{KL}%
\begin{tabular}[t]{l}
$=\{{\bf v}\in {\bf H}^{1}(\Omega );\,{\bf v}={\bf 0}$ on $%
\Gamma_{D},\,(s_{i3}({\bf v}))_{i=1,2,3}={\bf 0}\}$ \\
$=\{(\bar{v}_{1}-x_{3}{\partial }_{1}v_{3},\,\bar{v}_{2}-x_{3}{\partial }%
_{2}v_{3},\,v_{3});$
\begin{tabular}[t]{l}
${\bf \bar{v}:=(}\bar{v}_{1},\bar{v}_{2})\in {\bf H}^{1}({\omega }%
),\,v_{3}\in H^{2}({\omega }),$ \\
${\bf \bar{v}=\nabla }_{\hat{x}}v_{3}={\bf 0}\text{ and }v_{3}=0\;\text{on }%
\Gamma _{D}\}.$%
\end{tabular}%
\,%
\end{tabular}
\label{defi KL}
\end{equation}

\smallskip \noindent The following lemma was proved in [7], the
mean-operator ${\cal M}$ being defined in (\ref{1003}).

\medskip \noindent {\bf Lemma 4.1. }{\it If assumptions }(\ref{2}, \ref{8}, %
\ref{9}, \ref{49}) {\it and conventions} (\ref{1001}, \ref{492}, \ref{493})
{\it hold then}\vspace{0.15cm} \newline
{\rm {\it \smallskip \noindent }}{\it (i)\ for each fixed} $b$ {\it there is
a unique solution to} (\ref{63}); \vspace{0.15cm} \newline
{\rm {\it \smallskip \noindent }}{\it (ii)\ }$\Vert {\bf K}^{a}({\bf u}%
^{b})\Vert _{{\bf L}^{2}({\Omega })}+\Vert {\bf L}^{a}(\varphi ^{b})\Vert _{%
{\bf L}^{2}({\Omega }_{1}^{\varepsilon })}+$ $\sqrt{G_{1}}$ $\Vert \nabla _{%
\hat{x}}^{\varepsilon }{\cal M}(L_{3}^{a}(\varphi ^{b}))\Vert _{{\bf L}^{2}({%
\Omega }_{1}^{\varepsilon })}$ $\leq C$; \vspace{0.15cm} \newline
{\rm {\it \smallskip \noindent }}{\it (iii)\ there exists} ${\bf M=(K,L})\in
$ $(L^{2}({\Omega }\times Y))^{7}\times (L^{2}({\Omega }\times Y_{1}))^{3}$
{\it such that} $({\bf M}^{a}{\bf (U}^{b}))${\it \ two-scale converges to} $%
{\bf M}$ {\it in} ${\bf L}^{2}({\Omega }\times Y)\times {\bf L}^{2}({\Omega }%
\times Y_{1})$; \vspace{0.15cm} \newline
{\rm {\it \smallskip \noindent }}{\it (iv)\ } {\it there exists} ${\bf u}\in
{\bf V}_{KL}$ {\it and} ${\bf u}^{1}{\bf =}(u_{1}^{1},u_{2}^{1},0)$ {\it with%
} $u_{1}^{1},\,u_{2}^{1}\in L^{2}({\Omega };H_{\sharp
}^{1}(Y)/{\mathbb R})$ {\it such that} \vspace{0.15cm} \newline
\begin{tabular}[t]{l}
\thinspace \thinspace \thinspace \thinspace \thinspace \thinspace\ $({\bf u}%
^{b})$ {\it converges weakly to} ${\bf u}$ {\it in} ${\bf H}^{1}({\Omega })$%
, $\medskip $\vspace{0.05cm} \\
\thinspace \thinspace \thinspace \thinspace \thinspace \thinspace $(\nabla _{%
\hat{x}}{\bf u}^{b})${\bf {\rm {\it \ }}} {\it two-scale converges to }$%
\nabla _{\hat{x}}{\bf u+\nabla }_{y}{\bf u}^{1}$ {\it in }${\bf L}%
^{2}(\Omega \times Y_{1})$,$\medskip $ \vspace{0.05cm} \\
\thinspace \thinspace \thinspace \thinspace \thinspace \thinspace $(\partial
_{3}{\bf u}^{b})$ {\it two-scale converges to }$\partial _{3}{\bf u}$ {\it %
in }${\bf L}^{2}(\Omega \times Y_{1})$;%
\end{tabular}
\vspace{0.15cm} \newline
{\rm {\it \smallskip \noindent }}{\it (v)\ } $(\varphi ^{b})$ {\it two-scale
converges to }$\varphi _{m}$ {\it in }$L^{2}(\Omega \times Y_{1})$; \vspace{%
0.15cm} \newline {\rm {\it \smallskip \noindent }}{\it (vi)\ } {\it
there exists }$\varphi ^{1}\in L^{2}(\Omega ;H^{1}(Y_{1})/{\mathbb
R})$ {\it such that}$\ ^{t}(L_{1},L_{2})={\nabla }_{y}\varphi ^{1}$,
\vspace{0.15cm} \newline {\rm {\it \smallskip \noindent }}{\it
(vii)\ } ${\cal M}(L_{3})$ {\it is
independent of} $y$ {\it and } {\it for Dirichlet conditions }${\cal M}%
(L_{3})=\varphi _{c}$; \vspace{0.15cm} \newline
{\rm {\it \smallskip \noindent }}{\it (viii)\ In the case of nonlocal mixed
conditions}, ${\cal M}(L_{3})\in H^{1}(\omega )$ {\it and }$(\nabla _{\hat{x}%
}^{\varepsilon }{\cal M}(L_{3}^{a}(\varphi ^{b})))$ {\it two-scale converges
to } $\nabla _{\hat{x}}{\cal M}(L_{3})$ {\it in }${\bf L}^{2}(\Omega \times
Y_{1})$.

\bigskip \noindent {\bf 5. Limit Models I - Effective composite plate models
when the thickness is small with respect to the size of the inclusions. }

\smallskip \noindent This section is devoted to the statement of the
two-dimensional composite plate models when the thickness of the plate is
much smaller than the size of the inclusions, that is when $\varepsilon ,$ $%
a $ and $a/\varepsilon $ tend together to zero. The models are derived by
coupling the homogenization method and the asymptotic method for plates.
Consequently, they present characteristics of both approaches. Since no
other situation occurs when $a$ is small with respect to $\varepsilon $, it
would be a posteriori equivalent to derive first the two-dimensional plate
model and then to apply the method of homogenization to obtain a homogeneous
two dimensional plate model.

\smallskip \noindent The results are summarized in Theorem 5.1 below, which
is given after introducing a few notations: the projections involved by the
plate approach, and the local problems and homogenized tensors involved by
the homogenization process. The proof is postponed to Section 7. \bigskip

\bigskip \noindent {\bf 5.1 Notations related to plate theory }

\medskip \noindent The stiffness - piezoelectric - permittivity coefficients
of the two dimensional plate model obtained by eliminating the transverse
components are as follows.%
\begin{equation}
\left\{
\begin{array}{l}
\Pi \text{ and }\Pi _{1}\text{ are respectively the projections from }(L^{2}(%
{\Omega }\times Y))^{7}\times (L^{2}({\Omega }\times Y_{1}))^{3}\medskip \\
\text{on its subspaces }\left\{ ({\bf 0}_{4},(K_{i3})_{i=1..3},{\bf 0}%
_{2},L_{3})\right\} \text{ and }\left\{ ({\bf 0}_{9},L_{3})\right\} ,\text{ }%
\Pi _{2}=\Pi -\Pi _{1},\medskip \\
{\bf T}_{{\cal N}}=-(\Pi {\cal R}\Pi )^{-1}\Pi {\cal R}\text{, }{\bf T}_{%
{\cal M}}=-(\Pi _{2}{\cal R}\Pi _{2})^{-1}\Pi _{2}{\cal R}\,\text{,}\medskip
\\
{\cal R}_{{\cal N}}=(\text{Id}+{^{t}}{\bf T}_{{\cal N}}){\cal R}(\text{Id}+%
{\bf T}_{{\cal N}}),\text{ }{\cal R}_{{\cal M}}=(\text{Id}+{^{t}}{\bf T}_{%
{\cal M}})({\cal R}+2G\Pi _{1})(\text{Id}+{\bf T}_{{\cal M}}).%
\end{array}%
\right.  \label{Tenseurs}
\end{equation}%
{\bf Remark:} with notations like $(\Pi {\cal R}\Pi )^{-1}$, we mean the
inverse application of $\Pi {\cal R}\Pi $ as an application onto $\left\{ (%
{\bf 0}_{4},(K_{i3})_{i=1..3},{\bf 0}_{2},L_{3})\right\} $ (for the sake of
simplicity in the notation, $\Pi ,\Pi _{2}$ are identified with their
transposed applications).

\bigskip \noindent {\bf 5.2 Notations related to homogenization theory }

\smallskip \noindent Let $R_{{\cal M}{\alpha }\beta {\gamma }{\delta }},$
and $R_{{\cal N}{\alpha }\beta {\gamma }{\delta }}$ denote the relevant
coefficients of the tensors ${\cal R}_{{\cal M}},{\cal R}_{{\cal N}}$ and $%
{\cal R}_{{\cal M}}^{Mix}$ written on the format (\ref{tenseur}). Let
\[
H_{\sharp }^{2}(Y)=\{v\in H^{2}(Y);\text{ }v\text{ and }\nabla v\text{ \ are
}Y-\text{periodic}\}.
\]%
For any ${\bf v}\in {\bf H}^{1}({Y})$, let
\begin{equation}
S_{\alpha \beta }({\bf v)=}{\frac{1}{2}}({\partial }_{y_{\alpha }}v_{\beta }+%
{\partial }_{y_{\beta }}v_{\alpha }).  \label{Definition S}
\end{equation}%
Let $({\bf u}_{{\cal M}}^{{\gamma }{\delta }},u_{{\cal N}3}^{{\gamma }{%
\delta }})$ $\in (H_{\sharp }^{1}(Y))^{2}\times H_{\sharp }^{2}(Y),$ for $%
\gamma ,\delta \in \{1,2\},$ be the solutions of

\begin{equation}
\left\{
\begin{array}{l}
\forall {\bf v}\in (H_{\sharp }^{1}(Y))^{2},\ \ \dint\nolimits_{Y}S_{{\alpha
}\beta }({\bf v})R_{{\cal M}{\alpha }\beta {\lambda }\mu }S_{{\lambda }\mu }(%
{\bf u}_{{\cal M}}^{{\gamma }{\delta }})\,\,\text{d}y=-\dint\nolimits_{Y}%
\,S_{{\alpha }\beta }({\bf v})R_{{\cal M}{\alpha }\beta {\gamma }{\delta }}\,%
\text{d}y,\medskip \\
\forall v_{3}\in H_{\sharp }^{2}(Y),\ \ \dint\nolimits_{Y}{\partial }%
_{y_{\alpha }y_{\beta }}^{2}v_{3}R_{{\cal N}{\alpha }\beta {\lambda }\mu }{%
\partial }_{y_{{\lambda }}y_{\mu }}^{2}u_{{\cal N}3}^{{\gamma }{\delta }}\,\,%
\text{d}y=-\dint\nolimits_{Y}{\partial }_{y_{\alpha }y_{\beta }}^{2}v_{3}R_{%
{\cal N}{\alpha }\beta {\gamma }{\delta }}\,\text{d}y{\bf .}\medskip%
\end{array}%
\right.  \label{24}
\end{equation}%
Let ${\bf u}_{{\cal M}}^{3}$ $\in (H_{\sharp }^{1}(Y))^{2}$ be the solution
of

\begin{equation}
\forall {\bf v}\in (H_{\sharp }^{1}(Y))^{2},\text{ \ }\dint\nolimits_{Y}S_{{%
\alpha }\beta }({\bf v})R_{{\cal M}{\alpha }\beta {\lambda }\mu }S_{{\lambda
}\mu }({\bf u}_{{\cal M}}^{3})\,\text{d}y=-\dint\nolimits_{Y}d_{{\cal M}3{%
\alpha \beta }}S_{{\alpha }\beta }({\bf v})\,\text{d}y.  \label{241}
\end{equation}

\smallskip \noindent The effective tensors $R_{{\cal M}}^{H}$, $R_{{\cal N}%
}^{H}$, $d_{{\cal M}}^{H}$, $e_{{\cal M}}^{H}$ and $c_{33}^{H}$ are then
defined by%
\begin{equation}
\begin{array}{l}
R_{{\cal N}\gamma \delta \rho \xi }^{H}=\dint\nolimits_{Y}(\delta _{{\alpha }%
\beta ,\gamma \delta }+{\partial }_{y_{\alpha }y_{\beta }}^{2}u_{{\cal N}}^{{%
\gamma }{\delta }})R_{{\cal N}{\alpha }\beta {\lambda }\mu }(\delta _{{%
\lambda }\mu ,\rho \xi }+{\partial }_{y_{{\lambda }}y_{\mu }}^{2}u_{{\cal N}%
}^{\rho \xi })\text{ d}y,%
\end{array}
\label{TenseurRN}
\end{equation}%
\begin{equation}
\left\{
\begin{array}{l}
\left(
\begin{array}{cc}
R_{{\cal M}{\gamma }{\delta \rho \xi }}^{H} & d_{{\cal M}3{\gamma }{\delta }%
}^{H} \\
e_{{\cal M}3{\rho \xi }}^{H} & c_{{\cal M}33}^{H}%
\end{array}%
\right) =\dint_{Y}\left(
\begin{array}{cc}
\delta _{{\alpha }\beta ,\gamma \delta }+S_{{\alpha \beta }}({\bf u}_{{\cal M%
}}^{{\gamma }{\delta }}) & 0 \\
S_{{\lambda }\mu }({\bf u}_{{\cal M}}^{3}) & 1%
\end{array}%
\right) \medskip \\
\text{\ \ \ \ \ \ \ \ \ \ \ \ \ \ \ }\left(
\begin{array}{cc}
R_{{\cal M}{\alpha }\beta {\lambda }\mu } & d_{{\cal M}3{\alpha }\beta } \\
-d_{{\cal M}3{\lambda }\mu } & c_{{\cal M}33}%
\end{array}%
\right) \left(
\begin{array}{cc}
\delta _{{\lambda }\mu ,\rho \xi }+S_{{\lambda }\mu }({\bf u}_{{\cal M}%
}^{\rho \xi }) & S_{{\lambda }\mu }({\bf u}_{{\cal M}}^{3}) \\
0 & 1%
\end{array}%
\right) \text{ d}y{\bf ,}%
\end{array}%
\right.  \label{TenseursM}
\end{equation}%
where the notations $\delta _{{\lambda }\mu ,\rho \xi }$ are Kronecker
symbols, that is
\[
\delta _{{\lambda }\mu ,\rho \xi }%
\begin{tabular}[t]{l}
$=1\ $if $\lambda =\rho $ and $\mu =\xi ,$ \\
$=0\ $if not.%
\end{tabular}%
\]%
For local mixed conditions we also need%
\[
R_{{\cal M}{\alpha }\beta {\gamma }{\delta }}^{H,loc}=R_{{\cal M}{\alpha }%
\beta {\gamma }{\delta }}^{H}-\left( c_{{\cal M}33}^{H}+2|Y_{1}|G\right)
^{-1}d_{{\cal M}3{\alpha }\beta }^{H}e_{{\cal M}3{\gamma }{\delta }}^{H}.
\]%
\smallskip \noindent Let then
\begin{equation}
\ell _{u}({\bf v})=\dint\nolimits_{\Omega }{\bf f.v}\text{ d}%
x+\dint\nolimits_{\Gamma _{N}}{\bf g.v}\text{ d}s,  \label{ell u}
\end{equation}%
\[
\ell ({\bf v)}=\left\{
\begin{tabular}{l}
$\ell _{u}({\bf v})-2\dint_{\omega }s_{\alpha \beta }({\bf \bar{v}})d_{{\cal %
M}3\alpha \beta }^{H}\varphi _{c}$ d$\hat{x}$ \ (Dirichlet cond.), \\
$\ell _{u}({\bf v})-2\dint\nolimits_{{\omega }}s_{\alpha \beta }({\bf \bar{v}%
})d_{{\cal M}3\alpha \beta }^{H}$ $\left( c_{{\cal M}33}^{H}+2|Y_{1}|G%
\right) ^{-1}|Y_{1}|h$ d$\hat{x}$ \ (local mixed cond.), \\
$\ell _{u}({\bf v})+2\left\vert Y_{1}\right\vert \dint\nolimits_{{\omega }}%
\tilde{L}_{3}$ $h$ d$\hat{x}$\ (nonlocal mixed cond.).%
\end{tabular}%
\right.
\]

\bigskip \noindent {\bf 5.3 Theorem 5.1: effective models }

\medskip \noindent The set ${\bf V}_{KL}$ of Kirchoff-Love displacement
fields have been defined in (\ref{defi KL}).

\medskip \noindent {\bf Theorem 5.1{\rm {\it . }}}{\it Assume that the
assumptions of Section 2.4 hold and that }$a${\it , }$\varepsilon ${\it \
and }$a/\varepsilon ${\it \ tend together to zero,}

{\bf {\rm {\it \smallskip \noindent }}}({\it i) in the case of Dirichlet
conditions }$({\bf u}^{b})$ {\it converges to }%
\[
{\bf u}=(\bar{u}_{1}-x_{3}{\partial }_{1}u_{3},\bar{u}_{2}-x_{3}{\partial }%
_{2}u_{3},u_{3})\in {\bf V}_{KL}
\]%
{\it which is the unique solution in}{\bf {\rm {\it \ }}}${\bf V}_{KL}${\bf
{\rm {\it \ }}}{\it of}
\begin{equation}
\forall \text{ }{\bf v}\in {\bf V}_{KL},\ \ \dint\nolimits_{\omega }\left(
2s_{{\alpha }\beta }({\bf \bar{v}})R_{{\cal M}{\alpha }\beta {\gamma }{%
\delta }}^{H}s_{\gamma \delta }({\bf \bar{u}})+{\dfrac{2}{3}}\partial _{{%
\alpha }\beta }^{2}v_{3}R_{{\cal N}{\alpha }\beta {\gamma }{\delta }}^{H}{%
\partial }_{\gamma \delta }^{2}u_{3}\right) \text{ d}\hat{x}=\ell ({\bf v});
\label{EffectifDirichlet}
\end{equation}%
{\bf {\rm {\it \smallskip \noindent }}}({\it ii) in the case of local mixed
conditions, }$({\bf u}^{b})$ {\it converges to }${\bf u}$ {\it which is the
unique solution in }$V_{KL}${\it \ of}
\[
\forall \text{ }{\bf v}\in {\bf V}_{KL},\ \ \dint\nolimits_{{\omega }}\left(
2s_{{\alpha }\beta }({\bf \bar{v}})\ R_{{\cal M}{\alpha }\beta {\gamma }{%
\delta }}^{H,loc}\ s_{{\gamma \delta }}({\bf \bar{u}})+{\dfrac{2}{3}{%
\partial }_{{\alpha }\beta }^{2}}v_{3}R_{{\cal N}{\alpha }\beta {\gamma
\delta }}^{H}{\partial }_{{\gamma \delta }}^{2}u_{3}\right) \text{ d}\hat{x}%
=\ell ({\bf v});
\]%
{\bf {\rm {\it \smallskip \noindent }}}{\it (iii)\ for nonlocal mixed
conditions, }$({\bf u}^{b},{\cal M}(L_{3}^{a}(\varphi ^{b})))$ {\it %
converges to}
\[
({\bf u,}L_{3}^{0})=(\bar{u}_{1}-x_{3}{\partial }_{1}u_{3},\bar{u}_{2}-x_{3}{%
\partial }_{2}u_{3},u_{3},L_{3}^{0})\in {\bf V}_{KL}\times H^{1}(\omega )%
{\bf {\rm {\it \ }}}
\]%
{\it which is the unique solution of}%
\[
\left\{
\begin{array}{l}
\,\forall \,({\bf v,}\tilde{L}_{3})=(\bar{v}_{1}-x_{3}{\partial }_{1}v_{3},%
\bar{v}_{2}-x_{3}{\partial }_{2}v_{3},v_{3},\tilde{L}_{3})\in {\bf V}%
_{KL}\times H^{1}({\omega }),\medskip \\
\dint\nolimits_{{\omega }}\left( 2(s_{{\alpha }\beta }({\bf \bar{v}}),\tilde{%
L}_{3})\left(
\begin{array}{cc}
R_{{\cal M}{\alpha }\beta {\gamma }{\delta }}^{H} & d_{{\cal M}3{\alpha }%
\beta }^{H} \\
e_{{\cal M}3{\gamma }{\delta }}^{H} & c_{{\cal M}33}^{H}+2|Y_{1}|G%
\end{array}%
\right) \left(
\begin{array}{c}
s_{{\gamma \delta }}({\bf \bar{u}}) \\
L_{3}^{0}%
\end{array}%
\right) \right) \text{ d}\hat{x}\medskip \\
\ \ \ \ \ \ \ \ \ \ \ \ \ \ \ \ \ \ \ \ \ \ \ \ \ \ \ \ +\dint\nolimits_{{%
\omega }}\left( 4|Y_{1}|G_{1}{\partial }_{\alpha }\tilde{L}_{3}{\partial }%
_{\alpha }L_{3}^{0}+\dfrac{2}{3}{\partial }_{{\alpha }\beta }^{2}v_{3}R_{%
{\cal N}{\alpha }\beta {\gamma \delta }}^{H}{\partial }_{{\gamma \delta }%
}^{2}u_{3}\right) \text{ }\,\text{d}\hat{x}=\ell ({\bf v}).%
\end{array}%
\right.
\]

\noindent {\bf Remarks: }(i) Since elasticity coefficients are constant in
the thickness, the membrane and the flexion models are discoupled as for a
single-layer elastic plate.$\medskip $

\noindent (ii) for Dirichlet, and local mixed conditions, the limit model
has the same form as the elastic plate model, the influence of piezoelectric
inclusion only appears in the definition of the effective coefficients and
as a source term.$\medskip $

\noindent (iii) The piezoelectric force operates only in the membrane model
as for the single-layer piezoelectric plate. A piezoelectric force may be
generated in the flexion model by a torque in a multi-layer piezoelectric
plate.$\medskip $

\noindent (iv) In a thin plate with an imposed voltage $\varphi _{c}$ on the
upper face of a piezoelectric inclusion $\omega _{P}$, the resulting force
is concentrated along the boundary of the inclusion and may be expressed as
a distribution $v\mapsto \sum_{\alpha ,\beta =1}^{2}\left\langle d_{\alpha
\beta }\varphi _{c},s_{\alpha \beta }(v)\right\rangle _{\omega _{P}}$ where $%
v$ is a regular test function, $d$ the piezoelectric coefficient, and $%
\left\langle .,.\right\rangle _{\omega _{P}}$ the distribution bracket over
the mean surface of the inclusion. For a periodic distribution of
piezoelectric inclusions with an imposed voltage $\varphi _{c}$ on the upper
faces, as considered in this paper, the resulting force is distributed over
the complete mean surface $\omega $ and is still a distribution $v\mapsto
\sum_{\alpha ,\beta =1}^{2}\left\langle d_{\alpha \beta }\varphi
_{c},s_{\alpha \beta }(v)\right\rangle _{\omega }$ because in general, the
distributed voltage $\varphi _{c}\in L^{2}(\omega )$ has no regularity.$%
\medskip $

\noindent (v) For nonlocal conditions, the presence of derivatives in the
term $\dint\nolimits_{{\omega }}{\partial }_{\alpha }\tilde{L}_{3}{\partial }%
_{\alpha }L_{3}$ d$\hat{x}$ is due to the electric circuits connecting
neighbouring cells, which results in a regularization of the transverse
component of the electric field $L_{3}$ and in a direct nontrivial coupling
between mechanical and electrical effects. Note that one may a priori choose
the form of the differential operator on $L_{3}$ by choosing the way to
connect the upper faces of the inclusions to each others. This last point is
interresting in view of controlling the structure.$\medskip $

\noindent (vi) For the case of local mixed conditions, the local circuit,
i.e. without connections between cells, allows elimination of $L_{3}$ from
the effective equations. Even if it is actuated by the current source $h$,
the piezoelectric force is similar to the case of an imposed voltage with $%
\varphi _{c}=|Y_{1}|h/\left( c_{{\cal M}33}^{H}+2|Y_{1}|G\right) $ where $%
1/\left( c_{{\cal M}33}^{H}+2|Y_{1}|G\right) $ is the impedance of the local
circuit.$\medskip $

\noindent (vii) It would also be possible to keep the voltage $L_{3}$ in the
model with local mixed conditions. The resulting model would be similar to
the model with nonlocal condition, i.e. with the current source $h$, but
with $G_{1}=0$. This formulation is most suited for extension to dynamic
problems where the admittance $G$ is generally an integro-differential
operator in time.$\medskip $

\bigskip \noindent {\bf 6. Limit Models II - Effective composite plate
models when }$a\sim \varepsilon ${\bf . }

\smallskip \noindent This section is devoted to the statement of the
two-dimensional limit models when $a$ and $\varepsilon $ tend to zero while $%
a/\varepsilon $ tends to a finite positive limit. Then, there is no
loss of generality in assuming that $a/\varepsilon $ tends to $1$.
If $a/\varepsilon $ tends to any other positive number $\ell \in
{\mathbb R}^{+\ast }$, one just has to normalize the coordinates
($x_{1},x_{2}$) in an appropriate way. However $\ell $ must be {\em
not too far} from 1. If not, one of the two other models,
$a/\varepsilon \rightarrow 0$ if $\ell <<1$, $\varepsilon
/a\rightarrow 0$ if $\ell >>1$, would be to consider.

\smallskip \noindent The results are summarized in Theorem 6.1 below, which
is given after introducing necessary preliminary notations, in particular
the local problems to be solved to compute the effective coefficients of the
limit models. The proof is postponed to Section 8. \bigskip

\medskip \noindent {\bf 6.1 Notations related to the microstructure. \ }
Because $a$ and $\varepsilon $ are of the same order of magnitude, we have
here a unique three-dimensional microstructure containing both the $y$ and $%
x_{3}$ directions. Let us introduce appropriate notations%
\[
\begin{array}{l}
Z=Y\times ]-1,1[,\,Z_{1}=Y_{1}\times ]-1,1[,\,\,z=(y_{1},y_{2},x_{3}),%
\end{array}%
\]%
\[
\Gamma ^{+}=Y_{1}\times \left\{ 1\right\} ,\ \Gamma ^{-}=Y_{1}\times \left\{
-1\right\} ,
\]%
\begin{equation}
\begin{tabular}{l}
${\bf W}^{1}=\left( L^{2}(]-1,1[;H_{\sharp }^{1}(Y))\cap H^{1}(Z)\right)
^{3}\times \Psi ^{1}$ where$\medskip $ \\
$\Psi ^{1}=\left\{ \psi \in H^{1}(Z_{1});\ \psi =0\ \text{on }\Gamma
^{-}\cup \Gamma ^{+}\right\} $.%
\end{tabular}
\label{def espaces locaux}
\end{equation}%
\begin{equation}
\begin{array}{l}
\forall {\bf V}^{1}=({\bf v}^{1},\psi ^{1})\in {\bf W}^{1},\text{ }{\bf M}%
^{1}{\bf (V)}=(\left( S_{ij}({\bf v}^{1})\right) _{i,j=1..3},\nabla _{z}\psi
^{1}))\medskip \\
\text{where }S_{ij}({\bf v}^{1})={\dfrac{1}{2}}({\partial }_{z_{i}}v_{j}^{1}+%
{\partial }_{z_{j}}v_{i}^{1}).%
\end{array}
\label{Defi M1 et Sz}
\end{equation}%
\medskip \noindent Let ${\bf U}_{{\cal M}}^{\lambda \mu }=({\bf u}_{{\cal M}%
}^{\lambda \mu },{\varphi }_{{\cal M}}^{\lambda \mu })\in {\bf W}^{1}$
designates for any $\left( \lambda ,\mu \right) \in \left\{ 1,2\right\} ^{2}$
the solutions to%
\begin{equation}
\forall {\bf V}^{1}\in {\bf W}^{1},\ \dint\nolimits_{Z}{\bf M}^{1}{\bf (V}%
^{1}){\cal R\ }{^{t}}{\bf M}^{1}{\bf (U}_{{\cal M}}^{\lambda \mu })\text{ d}%
z=-\dint\nolimits_{Z}{\bf M}^{1}{\bf (V}^{1})\ ^{t}{\bf X}_{\lambda \mu }%
\text{ d}z,  \label{94}
\end{equation}%
where ${\bf X}_{\lambda \mu }=\left( R_{{\alpha }\beta \lambda \mu },2R_{{%
\alpha }3\lambda \mu },R_{33\lambda \mu },-d_{{\alpha }\lambda \mu
},-d_{3\lambda \mu }\right) $ is for the first column of ${\cal R}$.

\medskip \noindent Similarly, let ${\bf U}_{{\cal N}}^{\lambda \mu }=({\bf u}%
_{{\cal N}}^{\lambda \mu },{\varphi }_{{\cal N}}^{\lambda \mu })\in {\bf W}%
^{1}$ be solutions to%
\begin{equation}
\forall {\bf V}^{1}\in {\bf W}^{1},\ \dint_{Z}{\bf M}^{1}{\bf (V}^{1}){\cal %
R\ }{^{t}}{\bf M}^{1}{\bf (U}_{{\cal N}}^{\lambda \mu })\text{ d}%
z=\dint\nolimits_{Z}\text{ }x_{3}\text{ }{\bf M}^{1}{\bf (V}^{1})\ ^{t}{\bf X%
}_{\lambda \mu }\text{ d}z,  \label{95}
\end{equation}%
and ${\bf U}^{3}=({\bf u}^{3},{\varphi }^{3})\in {\bf W}^{1}$ be solution to%
\begin{equation}
\forall {\bf V}^{1}\in {\bf W}^{1},\ \dint_{Z}{\bf M}^{1}{\bf (V}^{1}){\cal R%
}\text{\ }{^{t}}{\bf M}^{1}{\bf (U}^{3})\text{ d}z=\dint_{Z}\text{ }{\bf M}%
^{1}{\bf (V}^{1})\ \ ^{t}{\bf X}_{3}\text{ d}z\text{ },\medskip  \label{96}
\end{equation}%
where ${\bf X}_{3}=\left( d_{3{\alpha }\beta },2d_{3{\alpha }3},d_{333},c_{{%
\alpha }3},c_{33}\right) $ stands for the last column of ${\cal R}$.

\smallskip \noindent The effective coefficients are then given by:

\[
\left\{
\begin{array}{l}
R_{{\cal MM}\lambda \mu \rho \xi }^{H}=\dint_{Z}{(}{\bf E}^{\lambda \mu }+%
{\bf M}^{1}{\bf (U}_{{\cal M}}^{\lambda \mu })){\cal R}\text{ }{^{t}}({\bf E}%
^{\rho \xi }+{\bf M}^{1}{\bf (U}_{{\cal M}}^{\rho \xi }))\text{ d}z,\medskip
\\
R_{{\cal MN}\lambda \mu \rho \xi }^{H}=\dint_{Z}{(}{\bf E}^{\lambda \mu }+%
{\bf M}^{1}{\bf (U}_{{\cal M}}^{\lambda \mu })){\cal R}\text{ }{^{t}}(-x_{3}%
{\bf E}^{\rho \xi }+{\bf M}^{1}{\bf (U}_{{\cal N}}^{\rho \xi }))\text{ d}%
z,\medskip \\
R_{{\cal NM}\lambda \mu \rho \xi }^{H}=\dint_{Z}{(}-x_{3}{\bf E}^{\lambda
\mu }+{\bf M}^{1}{\bf (U}_{{\cal N}}^{\lambda \mu })){\cal R}\text{ }{^{t}}(%
{\bf E}^{\rho \xi }+{\bf M}^{1}{\bf (U}_{{\cal M}}^{\rho \xi }))\text{ d}%
z,\medskip \\
R_{{\cal NN}\lambda \mu \rho \xi }^{H}=\dint_{Z}{(}-x_{3}{\bf E}^{\lambda
\mu }+{\bf M}^{1}{\bf (U}_{{\cal N}}^{\lambda \mu })){\cal R}\text{ }{^{t}}%
(-x_{3}{\bf E}^{\rho \xi }+{\bf M}^{1}{\bf (U}_{{\cal N}}^{\rho \xi }))\text{
d}z,%
\end{array}%
\right.
\]%
\[
\left\{
\begin{array}{l}
d_{{\cal MM}3\lambda \mu }^{H}=\dint_{Z}\left( {\bf E}^{\lambda \mu }+{\bf M}%
^{1}{\bf (U}_{{\cal M}}^{\lambda \mu })\right) {\cal R}\text{ }{^{t}}\left(
{\bf b}+{\bf M}^{1}{\bf (U}^{3})\right) \text{ d}z,\medskip \\
d_{{\cal NM}3\lambda \mu }^{H}=\dint_{Z}\left( -x_{3}{\bf E}^{\lambda \mu }+%
{\bf M}^{1}{\bf (U}_{{\cal N}}^{\lambda \mu })\right) {\cal R}\text{ }{^{t}}%
\left( {\bf b}+{\bf M}^{1}{\bf (U}^{3})\right) \text{ d}z,\medskip \\
e_{{\cal MM}3\alpha \beta }^{H}=\dint_{Z}\left( {\bf b}+{\bf M}^{1}{\bf (U}%
^{3})\right) {\cal R}\text{ }{^{t}}({\bf E}^{\alpha \beta }+{\bf M}^{1}{\bf %
(U}_{{\cal M}}^{\alpha \beta }))\text{ d}z,\medskip \\
e_{{\cal MN}3\alpha \beta }^{H}=\dint_{Z}\left( {\bf b}+{\bf M}^{1}{\bf (U}%
^{3})\right) {\cal R}\text{ }{^{t}}\left( -x_{3}{\bf E}^{\alpha \beta }+{\bf %
M}^{1}{\bf (U}_{{\cal N}}^{\alpha \beta })\right) \text{ d}z,\medskip \\
c_{{\cal MM}33}^{H}=\dint_{Z}\left( {\bf b}+{\bf M}^{1}{\bf (U}^{3})\right)
{\cal R}\text{ }{^{t}}\left( {\bf b}+{\bf M}^{1}{\bf (U}^{3})\right) \text{ d%
}z.%
\end{array}%
\right.
\]%
where ${\bf E}^{\lambda \mu }=\left( (\delta _{\alpha \beta ,\lambda \mu })_{%
{\alpha },\beta =1,2},{\bf 0}_{6}\right) $ and ${\bf b=}{(}{\bf 0}_{9},1)$%
.\bigskip

\bigskip \noindent {\bf 6.2 Effective models: Theorem 6.1}

\medskip \noindent The set ${\bf V}_{KL}$ of Kirchoff-Love displacement
fields have been defined in (\ref{defi KL}).

\medskip \noindent {\bf Theorem 6.1.{\rm {\it \ }}}{\it Assume that the
assumptions of Section 2.4 hold and in addition that }$a,\varepsilon $ {\it %
tend to zero and} $a/\varepsilon $ {\it tends to} $1${\it \ then}

{\it \smallskip \noindent (i) in the case of Dirichlet conditions, the limit}%
{\bf {\rm {\it \ }}}${\bf u}$ {\it of }$({\bf u}^{b})_{b>0}$ {\it satisfies:
}${\bf u=}(\bar{u}_{1}-x_{3}\partial _{1}u_{3}$ $,\bar{u}_{2}-x_{3}\partial
_{2}u_{3},u_{3})\in {\bf V}_{KL}${\bf {\rm {\it \ }}}{\it and is the unique
solution in }${\bf V}_{KL}$ {\it of}

\[
\left\{
\begin{array}{l}
\forall {\bf v}\in V_{KL}\text{,\ \ \ }\dint\nolimits_{\omega }(s_{\alpha
\beta }({\bf \bar{v})},\partial _{\alpha \beta }^{2}v_{3})\left(
\begin{array}{cc}
R_{{\cal MM}\alpha \beta \gamma \delta }^{H} & R_{{\cal MN}\alpha \beta
\gamma \delta }^{H} \\
R_{{\cal NM}\alpha \beta \gamma \delta }^{H} & R_{{\cal NN}\alpha \beta
\gamma \delta }^{H}%
\end{array}%
\right) \left(
\begin{array}{c}
s_{\gamma \delta }({\bf \bar{u})} \\
\partial _{\gamma \delta }^{2}u_{3}%
\end{array}%
\right) \text{ d}\hat{x}=\medskip \\
\ \ \ \ \ \ \ \ \ \ \ \ \ \ \ \ \ \ \ \ \ \ \ \ \ \ \ \ \ \ \ \ \ \ \ \ \ \
\ \ \ \ \ \ \ \ell _{u}({\bf v})-\dint\nolimits_{\omega }(s_{\alpha \beta }(%
\overline{{\bf v}}{\bf )},\partial _{\alpha \beta }^{2}v_{3})d_{{\cal MM}%
3\alpha \beta }^{H}\varphi _{c}\text{ d}\hat{x};%
\end{array}%
\right.
\]

\smallskip \noindent {\it (ii) in the case of mixed conditions, the limit}%
{\bf {\rm {\it \ }}}${\bf u=}(\bar{u}_{1}-x_{3}\partial _{1}u_{3}$ $,\bar{u}%
_{2}-x_{3}\partial _{2}u_{3},u_{3})$ {\it of }$({\bf u}^{b})_{b>0}$ {\it %
belongs to }${\bf V}_{KL}$, {\it the limit} $L_{3}^{0}$ {\it of} ${\cal M}%
(L_{3}(\varphi ^{b}))$ {\it belongs to} $L^{2}(\omega )${\it \ or to }$%
H^{1}(\omega )${\it \ for nonlocal conditions, and} $\left( {\bf u}%
,L_{3}^{0}\right) $ {\it is the unique solution in }${\bf V}_{KL}\times
L^{2}(\omega )$ {\it or in} ${\bf V}_{KL}\times H^{1}(\omega )$ {\it for
nonlocal conditions of}

\[
\left\{
\begin{array}{l}
\forall {\bf (v,}\widetilde{L}_{3})\in {\bf V}_{KL}\times L^{2}(\omega )(%
\text{{\it resp. }}{\bf V}_{KL}\times H^{1}(\omega )),\medskip \\
\dint\nolimits_{\omega }(s_{\alpha \beta }({\bf \bar{v})},\widetilde{L}%
_{3},\partial _{\alpha \beta }^{2}v_{3})\left(
\begin{array}{ccc}
R_{{\cal MM}\alpha \beta \gamma \delta }^{H} & d_{{\cal MM}3\alpha \beta
}^{H} & R_{{\cal MN}\alpha \beta \gamma \delta }^{H} \\
e_{{\cal MM}3\alpha \beta }^{H} & c_{{\cal MM}33}^{H} & d_{{\cal MN}3\alpha
\beta }^{H} \\
R_{{\cal NM}\alpha \beta \gamma \delta }^{H} & e_{{\cal MN}3\alpha \beta
}^{H} & R_{{\cal NN}\alpha \beta \gamma \delta }^{H}%
\end{array}%
\right) \left(
\begin{array}{c}
s_{\gamma \delta }({\bf \bar{u})} \\
L_{3}^{0} \\
\partial _{\gamma \delta }^{2}u_{3}%
\end{array}%
\right) \text{ d}\hat{x}\medskip \\
\ \ \ \ \ \ \ \ \ \ \ \ \ \ \ +\dint\nolimits_{\omega }\left( 4|Y_{1}|(G%
\widetilde{L}_{3}L_{3}^{0}\,+G_{1}{\partial }_{\alpha }\widetilde{L}_{3}{%
\partial }_{\alpha }L_{3}^{0})\right) \text{d}\hat{x}=\ell _{u}({\bf v}%
)+2|Y_{1}|\dint_{{\omega }}\widetilde{L}_{3}h\,\text{d}\hat{x}\text{. }%
\end{array}%
\right.
\]

\noindent {\bf Remarks: }(i) On the contrary to the model in Section 5, the
membrane and flexion models are coupled. This comes from the presence of a
microstructure which size is comparable to the plate thickness, and
generates a complex displacement field at the scale of the plate thickness.

\noindent (ii) The other remarks (iii-vii) also hold for this model.

\bigskip

\bigskip \noindent {\bf 7. Proof of Theorem 5.1. }

\smallskip \noindent This section is devoted to the derivation of Theorem
5.1. The proof is based on the general results of Lemma 4.1. The proof
consists in three steps.\ The first one is the characterization of the limit
${\bf M}$ of $({\bf M}^{a}({\bf V))}$ and is the aim of Section 7.1. The
second one consists in the elimination of the $x_{3}-$variable as usual in
plate theory. The last one is the elimination of the local variable as in
homogenization theory and obtention of the effective models. The last two
steps slightly differ from one boundary condition to the other, and thus are
presented separately. In particular, in the case of local mixed conditions,
the complete elimination of the transverse component is achieved after the
elimination of the local variable. This two steps are presented in Section
7.2. The assumption that ${\cal R}$ does not depend on $x_{3}$ is widely
used there to obtain simplifications and compute our limit models.

\bigskip \noindent {\bf 7.1 Step 1: characterization of the limit }${\bf M}$%
{\bf . }

\smallskip \noindent In the forthcoming Lemma 7.4, the limit ${\bf M}$ is
related to the limits ${\bf u}$, ${\bf u}^{1}$ and $\varphi ^{1}$ defined in
Lemma 4.1. The limit of equation (\ref{63}) is given in Lemma 7.5. Suitable
choices of tests functions ${\bf V}$ lead here to a list of equations
corresponding to the various asymptotic levels of ${\bf M}^{a}({\bf V).}$
But first, we give \`{a} few additionnal notations, and technical results,
in section 7.1.1.

\bigskip \noindent {\bf 7.1.1 Further notations, and preliminary lemmas.}

\smallskip \noindent Recall that $C_{\sharp }^{\infty }(Y)$ denotes the
subspace of all $Y-$periodic functions of $C^{\infty }({\mathbb
R}^{2})$.
Similarly, let $C_{\sharp }^{\infty }(Y_{1})$ denote the subspace of all $Y-$%
periodic functions of $C^{\infty }({\mathbb Z}^{2}+Y_{1})$. For
functions $v$ on $\bar{\Omega}\times {\mathbb R}^{3}$ which are
$Y-$periodic with respect to
the second variable, we designate by $v^{\varepsilon }$ the function $%
x\longmapsto v(x,\hat{x}/\varepsilon )$. We also use test functions in%
\begin{equation}
{\bf W}_{ad}^{1}=\left\{ ({\bf v}^{1},\psi ^{1})\in {\bf D}(\bar{\Omega}%
,C_{\sharp }^{\infty }(Y))\times \Psi _{ad}^{1};\ {\bf v}^{1}={\bf 0}\text{
on }\Gamma _{D}\times Y\right\} ,  \label{W1ad}
\end{equation}%
\[
{\bf W}_{ad}=\left\{ ({\bf v},\psi )\in {\bf H}^{1}(\Omega )\times \Psi _{ad}%
\text{; }{\bf v}={\bf 0}\text{ on }\Gamma _{D}\right\} ,
\]%
where%
\[
\Psi _{ad}^{1}%
\begin{tabular}[t]{l}
$={\cal D}(\omega \times ]-1,1[;C_{\sharp }^{\infty }(Y_{1}))\ \ \ $for
Dirichlet conditions,$\medskip $ \\
$=\left\{ \psi ^{1}\in {\cal D}(\omega \times ]-1,1];C_{\sharp }^{\infty
}(Y_{1}));\ \psi ^{1}\text{ is constant for }x_{3}=1\right\} \ \ \ $%
otherwise,%
\end{tabular}%
\]%
\[
\Psi _{ad}%
\begin{tabular}[t]{l}
$={\cal D}\left( \omega \times ]-1,1[\right) \ \ \ \text{for Dirichlet
conditions,}\medskip $ \\
$=\left\{ \psi \in {\cal D}\left( \omega \times ]-1,1]\right) ;\ \psi \text{
is constant for }x_{3}=1\right\} \ \ \ $otherwise.%
\end{tabular}%
\]%
Functions in ${\bf W}_{ad}$ and in ${\bf W}_{ad}^{1}$ are admissible test
functions for (\ref{63}).

\noindent For ${\bf V=}({\bf v},\psi )\in {\bf W}_{ad}^{1}$, we have the
natural decomposition of ${\bf M}^{a}({\bf V}^{\varepsilon })$:
\begin{equation}
{\bf M}^{a}({\bf V}^{\varepsilon })=({\bf M}^{00}({\bf V}))^{\varepsilon }+%
\frac{1}{\varepsilon }({\bf M}^{10}({\bf V}))^{\varepsilon }+\frac{1}{a}(%
{\bf M}^{01}({\bf V}))^{\varepsilon }+\frac{1}{\varepsilon a}({\bf M}^{11}(%
{\bf V}))^{\varepsilon }+\frac{1}{a^{2}}({\bf M}^{02}({\bf V}))^{\varepsilon
},  \label{1331}
\end{equation}%
where%
\begin{equation}
\left\{
\begin{array}{l}
{\bf M}^{00}({\bf V})=((s_{\alpha \beta }({\bf v))}_{\alpha ,\beta =1,2},%
{\bf 0}_{3},(\partial _{\alpha }\psi )_{\alpha =1,2},0),\medskip \\
{\bf M}^{10}({\bf V})=((S_{\alpha \beta }({\bf v))}_{\alpha ,\beta =1,2},%
{\bf 0}_{3},(\partial _{y_{\alpha }}\psi )_{\alpha =1,2},0),\medskip \\
{\bf M}^{01}({\bf V})=({\bf 0}_{2\times 2},(s_{\alpha 3}({\bf v))}_{\alpha
=1,2},{\bf 0}_{3},\partial _{3}\psi ),\medskip \\
{\bf M}^{11}({\bf V})=({\bf 0}_{2\times 2},({\dfrac{1}{2}}\partial
_{y_{\alpha }}v_{3}{\bf )}_{\alpha =1,2},{\bf 0}_{4}),\medskip \\
{\bf M}^{02}({\bf V})=({\bf 0}_{2\times 2},,{\bf 0}_{2},s_{33}({\bf v)},{\bf %
0}_{3}).%
\end{array}%
\right.  \label{2331}
\end{equation}%
Notations $s_{ij}({\bf v)}$ and $S_{\alpha \beta }({\bf v)}$ have been
defined in (\ref{1005}) and (\ref{Definition S}). Relevant subspaces ${\mathbb M%
}$, ${\mathbb M}^{0}$, ${\mathbb M}^{-1}$ and ${\mathbb M}^{-2}$
associated to this decomposition are defined by
\begin{equation}
\left\{
\begin{array}{l}
{\mathbb M}^{0}=\{(K_{\alpha \beta },{\bf 0}_{6});\,\,K_{{\alpha }\beta }=s_{{%
\alpha }\beta }({\bf v})+S_{{\alpha }\beta }({\bf v}^{1}),\ {\bf v\in V}_{KL}%
{\bf ,}\text{ }{\bf v}^{1}\in {\bf V}_{KL}^{1}\},\medskip \\
{\mathbb M}^{-1}=\{({\bf 0}_{2\times 2},K_{{\alpha }3},{\bf 0}%
_{3},L_{3});\,\,K_{{\alpha }3}\in L^{2}(\Omega \times Y),\text{ }L_{3}\in
L^{2}({\Omega }\times Y_{1}),\text{ }\medskip \\
{\cal M}(L_{3})\text{ is independant of \ }y\text{; and }{\cal M}(L_{3})=0%
\text{ for Dirichlet conditions}\},\medskip \\
{\mathbb M}^{-2}=\{({\bf 0}_{2\times 2},{\bf 0}_{2},K_{33},{\bf 0}%
_{3});\,K_{33}\in L^{2}(\Omega \times Y)\},\medskip \\
{\mathbb M}={\mathbb M}^{-2}\oplus {\mathbb M}^{-1}\oplus {\mathbb M}^{0}.%
\end{array}%
\right.  \label{71}
\end{equation}%
where%
\begin{equation}
{\bf V}_{KL}^{1}=\{((\bar{v}_{\alpha }^{1}-x_{3}\partial _{y_{\alpha
}}v_{3}^{2})_{{\alpha }=1,2},0);\bar{v}_{\alpha }^{1}\in L^{2}(\omega
;H_{\sharp }^{1}(Y)),v_{3}^{2}\in L^{2}(\omega ;H_{\sharp }^{2}(Y))\}.
\label{defiKL1}
\end{equation}%
Our first lemma contains density results that will lead to a suitable weak
formulation for the problem solved by ${\bf M}$.

\bigskip \noindent {\bf Lemma 7.1. }{\it (i) The set }$\{{\bf M}^{02}({\bf V}%
);\,\,{\bf V}\in {\bf W}_{ad}^{1}\}${\it $\,\,\text{is dense in}\,$}$\,{\mathbb %
M}^{-2}${\it ,}

{\it \smallskip \noindent (ii) The set }$\{{\bf M}^{01}({\bf V});\,\,{\bf V}%
\in {\bf W}_{ad}^{{}},\,\,v_{3}=0\}\,\,\text{is dense in}\,\,{\mathbb M}^{-1}$,%
{\it \ }

{\it \smallskip \noindent (iii) The set }$\{{\bf M}^{0}({\bf V})+{\bf M}%
^{10}({\bf V}^{1});\,\,({\bf V},{\bf V}^{1})=(({\bf v,}0),({\bf v}^{1},0))$
{\it where} $({\bf V,V}^{1})\in $ ${\bf W}_{ad}^{1},({\bf v,v}^{1})$ $\in $ $%
{\bf V}_{KL}\times {\bf V}_{KL}^{1}\}$ {\it is dense in }${{\mathbb M}^{0}}$%
{\it .}

\bigskip \noindent {\bf Proof }:{\rm {\it {\bf \ }}}(i) Let $K_{33}\in {\cal %
C}^{\infty }(\Omega \times Y)$ with compact support. This function $K_{33}$
can also be considered as a function of ${\cal D}(\overline{\Omega },{\cal C}%
_{\sharp }^{\infty }(Y))$ by letting
\[
\forall x\in \Omega ,\ \forall y\in {\mathbb R}^{2}\text{, }%
K_{33}(x,y)=K_{33}(x,y^{\prime }),
\]%
$\ $where $y^{\prime }$ is d\'{e}fined by: $y^{\prime }\in Y$ and $%
y-y^{\prime }\in {\mathbb Z}^{2}$. Then,\ $v_{3}$ defined by:
\[
v_{3}(x,y)={\dint }_{0}^{x_{3}}K_{33}(\hat{x},t,y)\text{d}t
\]
satisfies: $({\bf 0}_{2},v_{3},0)\in {\bf W}_{ad}^{1}$ and $\partial
_{3}v_{3}=K_{33}$. As the set of functions of ${\cal C}^{\infty
}(\Omega \times Y)$ with compact support is dense in $L^{2}(\Omega
\times Y)$, this proves that $\{{\bf M}^{02}({\bf V});$ ${\bf V}\in
{\bf W}_{ad}^{1}\}$ is dense in ${\mathbb M}^{-2}$.

{\it \smallskip \noindent }(ii) The proof is similar for the first
two terms $K_{\alpha 3}$ of the elements ${\mathbb M}^{-1}$. For the
third term $L_{3}:$ let $L_{3}$ be a given function in ${\cal
C}^{\infty }(\Omega \times Y_{1})$ with compact support; we also
assume that ${\cal M}(L_{3})$ does not depend on $y$, and that in
addition ${\cal M}(L_{3})=0$ in the case of Dirichlet
conditions. As in (i), we now consider $L_{3}$ as a function of ${\cal D}(%
\overline{\Omega },{\cal C}_{\sharp }^{\infty }(Y_{1}))$; let then $\psi $
be defined by $\psi (x,y)={\dint }_{-1}^{x_{3}}L_{3}(\hat{x},t,y)$ d$t$.
With the two conditions above on $L_{3}$, $\psi $ is an admissible test
function such that $\partial _{3}\psi =L_{3}$; with densiry arguments again,
this completes the proof of (ii).

{\it \smallskip \noindent }Part (iii) is just a restatement of the
definition of ${{\mathbb M}^{0}.}$ \ \ \ $\Box $

\smallskip \noindent We now give two technical lemmas that are used in the
next section to caracterize the limit ${\bf M}$.

\bigskip \noindent {\bf Lemma 7.2.}{\it {\bf {\rm \ }}Let }$(u^{\varepsilon
})_{{\varepsilon }>0}${\it \ be bounded in }$H^{1}({\Omega })${\it . Let }$%
u\in H^{1}({\Omega })${\it \ and }$u^{1}\in L^{2}({\Omega };H_{\sharp
}^{1}(Y))${\it \ be functions such that }$(u^{\varepsilon })_{{\varepsilon }%
>0}${\it \ converges weakly to }$u${\it \ in }$H^{1}({\Omega })${\it , and }$%
(\nabla u^{\varepsilon })_{{\varepsilon }>0}${\it \ two-scale converges to }$%
\nabla u+\nabla _{y}u^{1}${\it . Then }

\[
\forall v\in {\it {\bf {\rm {\cal D}}}}(\Omega ;C_{\sharp }^{\infty }(Y)),\
\ \lim\limits_{\varepsilon \rightarrow 0}\int\nolimits_{{\Omega }}\frac{%
u^{\varepsilon }}{{\varepsilon }}\partial _{y_{\alpha }}v^{\varepsilon }\,%
\text{d}x=\int\nolimits_{{\Omega }\times Y}u^{1}\partial _{y_{\alpha }}v%
\text{ d}x\text{d}y{\rm .}
\]

\bigskip \noindent {\bf Proof :}{\it {\bf {\rm \ }}}One just has to pass to
the limit in
\[
\dint\nolimits_{{\Omega }}{\partial }_{\alpha }u^{\varepsilon }(x)v(x,%
\widehat{x}/{\varepsilon })\,\text{d}x=-\dint\nolimits_{{\Omega }%
}u^{\varepsilon }(x){\partial }_{\alpha }v(x,\widehat{x}/{\varepsilon })\,%
\text{d}x-{\dfrac{1}{{\varepsilon }}}\dint\nolimits_{{\Omega }%
}u^{\varepsilon }(x){\partial }_{y_{\alpha }}v(x,\widehat{x}/{\varepsilon }%
)\,\text{d}x.
\]

\medskip \noindent An integration by parts of the first term on the
right-hand side yields then the result. $\ \ \Box $

\medskip \noindent {\bf Lemma 7.3. }{\it The family}\ $\left( \dfrac{{%
\varepsilon }^{2}}{a}\dint\nolimits_{{\Omega }}\left( R_{\alpha
3kl}^{\varepsilon }K_{kl}^{a}({\bf u}^{b})+d_{\alpha 3k}^{\varepsilon
}L_{k}^{a}(\varphi ^{b})\right) \tilde{v}_{\alpha }\text{d}x\right) ${\it \
tends to zero for any} $\tilde{v}_{\alpha }\in {\cal D}(\bar{\omega}%
;C_{\sharp }^{\infty }(Y))$ {\it such that }$\tilde{v}_{\alpha }=0${\it \ on
}${\Gamma }_{D}${\it .}

\medskip \noindent {\bf Proof }{\rm {\it {\bf : }}}Choose\ ${\bf v}%
=\varepsilon ^{2}x_{3}$ $(\tilde{v}_{1},\tilde{v}_{2},0)$ and $\psi =0${\rm
{\it {\bf \ }}}in (\ref{63}), where{\rm {\it {\bf \ }}}$\tilde{v}_{\alpha
}\in {\cal D}(\bar{\omega};C_{\sharp }^{\infty }(Y))$ and $\tilde{v}_{\alpha
}=0$ on{\rm {\it {\bf \ }}}${\Gamma }_{D}$.{\rm {\it {\bf \ }}}Then
\[
\begin{array}{l}
\varepsilon ^{2}\dint\nolimits_{{\Omega }}(R_{\alpha \beta kl}^{\varepsilon
}K_{kl}^{a}({\bf u}^{b})+d_{\alpha \beta k}^{\varepsilon }L_{k}^{a}(\varphi
^{b}))(s_{{\alpha }\beta }({\bf v})+\dfrac{1}{\varepsilon }S_{{\alpha }\beta
}({\bf v}))\text{ d}x+\medskip \\
\ \ \ \ \ \ \ \ \ \ \ \ \ \ \ \ \ \ \ \varepsilon ^{2}a^{-1}\dint\nolimits_{{%
\Omega }}2(R_{\alpha 3kl}^{\varepsilon }K_{kl}^{a}({\bf u}^{b})+d_{\alpha
3k}^{\varepsilon }L_{k}^{a}(\varphi ^{b})){\rm {\it {\bf \ }}}\tilde{v}%
_{\alpha }\text{ d}x=\varepsilon ^{2}\ell _{u}^{b}({\bf v}).%
\end{array}%
\]

\smallskip \noindent Passing to the limit yields the result. $\ \ \ \Box $

\bigskip \noindent {\bf 7.1.2 Characterization of the limit }${\bf M}${\bf .
}

\medskip \noindent The first lemma of this section relates to the form of
the limit ${\bf M}$.\ We prove that ${\bf M\in }{\phi }_{c}+{\mathbb M}$ where $%
{\phi }_{c}=({\bf 0}_{9},\varphi _{c})$. The second gives a first
variational characterization of ${\bf M}$.

\medskip \noindent {\bf Lemma 7.4.} {\it Let }${\bf M}$ {\it designates the
two-scale limit, up to the extraction of a subsequence, of }$({\bf M}^{b}(%
{\bf U}^{b}))_{b>0}$. {\it Then }${\bf M\in }{\phi }_{c}+{\mathbb
M}${\it \
where} ${\phi }_{c}=({\bf 0}_{9},\varphi _{c})${\it , or in other words:}$%
\medskip $

{\it (i)\ \ there exist }${\bf u\in V}_{KL}$ and ${\bf u}^{1}\in {\bf V}%
_{KL}^{1}$ such that $K_{a\beta }=s_{\alpha \beta }({\bf u})+S_{\alpha \beta
}({\bf u}^{1});\medskip $

{\it (ii) }$L_{1}=L_{2}=0$;$\medskip $

{\it (iii) }$L_{3}\in L^{2}({\Omega }\times Y_{1}),$ ${\cal M}(L_{3})${\it \
is independant of }\ $y${\it ; and }${\cal M}(L_{3})=\varphi _{c}$ {\it for
Dirichlet conditions.}

\medskip \noindent The third point is just a restatement of Lemma 4.1 (vii).

\medskip \noindent {\bf Proof of (i)}: In view of Lemma 4.1, to prove (i),
it remains to prove that ${\bf u}^{1}\in {\bf V}_{KL}^{1}$. The choice $%
u_{3}^{1}=0$ has been done in Lemma 4.1 (iv). The computation of $u_{{\alpha
}}^{1}$ requires some calculus.

\smallskip \noindent First, with a few integration by parts we get
\begin{equation}
\text{%
\begin{tabular}{l}
$\dint\nolimits_{\Omega }(\partial _{{1}}u_{2}^{b}-\partial
_{2}u_{1}^{b})\partial _{3}v^{\varepsilon }\text{ d}x=$ \\
\ \ $2a\int\nolimits_{\Omega }K_{23}^{a}({\bf u}^{b})\left( \partial
_{1}v^{\varepsilon }+\dfrac{1}{\varepsilon }\partial _{y_{1}}v^{\varepsilon
}\right) \text{ d}x-2a\int\nolimits_{\Omega }K_{13}^{a}({\bf u}^{b})\left(
\partial _{2}v^{\varepsilon }+\dfrac{1}{\varepsilon }\partial
_{y_{2}}v^{\varepsilon }\right) \text{ d}x$%
\end{tabular}%
}  \label{eq u un}
\end{equation}%
for all $v\in {\cal D}(\Omega ,{\cal C}_{\sharp }^{\infty }(Y))$. Passing to
the limit, as $a/{\varepsilon }$ tends to zero, we get
\[
\forall v\in {\cal D}(\Omega ,{\cal C}_{\sharp }^{\infty }(Y)),\ \ \
\dint\nolimits_{\Omega \times Y}(\partial _{{1}}u_{2}-\partial
_{2}u_{1})\partial _{3}v\,\text{d}x\text{d}y+\dint\nolimits_{\Omega \times
Y}(\partial _{y_{1}}u_{2}^{1}-\partial _{y_{2}}u_{1}^{1})\partial _{3}v\,%
\text{d}x\text{d}y=0\text{.}
\]%
But because ${\bf u}\in {\bf V}_{KL}$, this reduces to%
\begin{equation}
\forall v\in {\cal D}(\Omega ,C_{\sharp }^{\infty }(Y)),\ \ \
\dint\nolimits_{\Omega \times Y}(\partial _{y1}u_{2}^{1}-\partial
_{y_{2}}u_{1}^{1})\partial _{3}v\,\text{d}x\text{d}y=0.  \label{26}
\end{equation}%
\smallskip \noindent Now, we show that\
\begin{equation}
\forall v\in {\cal D}(\Omega \times Y),\ \ \ \dint\nolimits_{\Omega \times
Y}\partial _{y_{\beta }}u_{\alpha }^{1}\partial _{33}^{2}v\,\text{d}x\text{d}%
y=0.  \label{27}
\end{equation}%
\smallskip \noindent As $(K_{33}^{a}({\bf u}^{b}))$ is bounded, $\varepsilon
^{-1}\dint\nolimits_{{\Omega }}{\partial }_{3}u_{3}^{b}({\partial }_{{\alpha
}}\tilde{v}^{\varepsilon }+\varepsilon ^{-1}{\partial }_{y_{\alpha }}\tilde{v%
}^{\varepsilon })\,$d$x$ tends to zero for any{\rm {\it {\bf \ }}}$\tilde{v}%
\in {\cal D}(\Omega ,{\cal C}_{\sharp }^{\infty }(Y))$; hence, $\varepsilon
^{-1}\dint\nolimits_{{\Omega }}{\partial }_{{\alpha }}u_{3}^{b}{\partial }%
_{3}\tilde{v}^{\varepsilon }\,$d$x$ tends to zero. Also, since $(K_{{\alpha }%
3}^{a})$ is bounded, $\varepsilon ^{-1}s_{{\alpha }3}({\bf u}^{b})$ tends to
zero. Consequently, $\varepsilon ^{-1}\dint\nolimits_{{\Omega }}{\partial }%
_{3}u_{{\alpha }}^{b}{\partial }_{3}\tilde{v}^{\varepsilon }\,$d$x$ tends to
zero. Choosing\ now $\tilde{v}={\partial }_{y_{\beta }}v$ where $v\in ${\rm
{\it {\bf ${\cal D}$}}}$(\Omega \times Y)${\rm {\it {\bf ,}}} one gets that $%
\varepsilon ^{-1}\dint\nolimits_{{\Omega }}u_{{\alpha }}^{b}{\partial }%
_{y_{\beta }}{\partial }_{33}^{2}v^{\varepsilon }\,$d$x$ tends to zero.
Then, using Lemma 7.2, we get (\ref{27}).

\medskip \noindent Now, we are able to conclude the proof. Equation (\ref{27}%
) implies that $\nabla _{y}u_{{\alpha }}^{1}$ is affine with respect to $%
x_{3}$. Thus, since $(u_{1}^{1},u_{2}^{1})$ is defined up to a function of $%
x $, $u_{{\alpha }}^{1}(x,y)$ = $\overline{u}_{{\alpha }}^{1}(\hat{x}%
,y)+x_{3}u_{{\alpha }}^{2}(\hat{x},y)$ where $\bar{u}_{{\alpha }}^{1},u_{{%
\alpha }}^{2}\in L^{2}({\omega }\times Y)$. But then ${\partial }%
_{3}u_{\alpha }^{1}=\overline{u}_{\alpha }^{1}\in L^{2}({\omega }\times Y)$.
On the other hand, in view of (\ref{26}), $(\partial _{3}u_{1}^{1},\partial
_{3}u_{2}^{1})$ is curl free, so that there exists $c_{\alpha }\in L^{2}({%
\omega })$ and $u_{3}^{2}\in L^{2}({\omega };H_{\sharp }^{1}(Y))$ such that $%
{\partial }_{3}u_{\alpha }^{1}$ = $c_{\alpha }-\partial _{y_{\alpha
}}u_{3}^{2}$. This finally implies that $(u_{1}^{1},u_{2}^{1})$ may be
chosen of the form $u_{{\alpha }}^{1}$ = $\bar{u}_{{\alpha }}^{1}-x_{3}{%
\partial }_{y_{\alpha }}u_{3}^{2}$ where ${\bf \bar{u}}^{1}$ $\in $ $(L^{2}({%
\omega }\times Y))^{2}$ and $u_{3}^{2}$ $\in L^{2}({\omega };H_{\sharp
}^{1}(Y))$. \ $\Box $

\medskip \noindent {\bf Proof of (ii)}: Passing to the limit in
\begin{equation}
a\int_{\Omega _{1}^{\varepsilon }}L_{3}^{a}(\varphi ^{b})\left( \partial
_{\alpha }\psi ^{\varepsilon }+\dfrac{1}{\varepsilon }\partial _{y_{\alpha
}}\psi ^{\varepsilon }\right) \text{ d}x=\int_{\Omega _{1}^{\varepsilon
}}L_{\alpha }^{a}(\varphi ^{b})\partial _{3}\psi ^{\varepsilon }\text{ d}x,
\label{Eq Lim L3}
\end{equation}%
because $a/\varepsilon $ tends to zero, we get that
\[
\forall {\psi }\in {\cal D}({\Omega }\times Y_{1})\text{, }\dint\nolimits_{{%
\Omega }\times Y_{1}}L_{\alpha }{\partial }_{3}\psi \text{ d}x\text{d}y=0.
\]%
This proves that $L_{1}$ and $L_{2}$ are independent of $x_{3}$. On the
other hand, a few integrations by parts, with test functions $\psi $
independent of $x_{3}$, yield%
\[
\begin{tabular}{l}
$a\dint_{{\Omega }_{1}^{\varepsilon }}L_{3}^{a}(\varphi ^{b})\,(x_{3}-1)\,({%
\partial }_{\alpha }{\psi }^{\varepsilon }+{\dfrac{1}{\varepsilon }}{%
\partial }_{y_{\alpha }}{\psi }^{\varepsilon })\,$d$x=$ \\
\ \ \ \ \ \ $\dint_{{\Omega }_{1}^{\varepsilon }}{\partial }_{\alpha
}\varphi ^{b}\psi ^{\varepsilon }\,$d$x+2\dint_{{\omega }_{1}^{\varepsilon
}}(\varphi _{m}^{b}-a\varphi _{c}^{b})({\partial }_{\alpha }{\psi }%
^{\varepsilon }+{\dfrac{1}{{\varepsilon }}}{\partial }_{y_{\alpha }}\psi
^{\varepsilon })\,$d$\hat{x}.$%
\end{tabular}%
\]%
The left-hand side tends to zero because $a,\varepsilon $ and $a/\varepsilon
$ tend to zero. Since $\varphi _{m}^{b}$ and $\varphi _{c}^{b}$ are constant
on each connected part of ${\omega }_{1}^{\varepsilon }$, the last term on
the right is zero. Therefore, passing to the limit, we get
\[
0=2\dint\nolimits_{{\omega }\times Y_{1}}L_{\alpha }{\psi }\text{ d}x\text{d}%
y.
\]

\smallskip \noindent Since $L_{\alpha }$ does not depend on $x_{3}$, this
proves that $L_{\alpha }=0$. $\ \Box $

\medskip \noindent {\bf Lemma 7.5 } {\it The family }$({\bf M}^{a}({\bf U}%
^{b}))_{b>0}$ {\it two-scale converges to} ${\bf M}$ {\it which is
the unique solution in }${\phi }_{c}+{\mathbb M}${\it \ of}

\begin{equation}
\left\{
\begin{array}{l}
\forall \,{\widetilde{{\bf M}}=(\widetilde{{\bf K}},\widetilde{{\bf L}})}\in
{\mathbb M},\ \ \ \dint_{{\Omega }\times Y}{\widetilde{{\bf M}}\ }{\cal R\ }{%
^{t}}{\bf M}\,\text{d}x\text{d}y+2G\dint_{{\Omega }\times Y_{1}}{\cal M}%
(L_{3}){\cal M}(\tilde{L}_{3})\,\text{d}x\text{d}y\medskip \\
\ \ \ \ \ \ \ \ \ \ \ \ \ +2G_{1}\dint_{{\Omega }\times Y_{1}}{\partial }%
_{\alpha }{\cal M}(L_{3}){\partial }_{\alpha }{\cal M}(\tilde{L}_{3})\,\text{%
d}x\text{d}y=\ell _{u}({\bf v})+\ell _{\varphi }(\tilde{L}_{3})%
\end{array}%
\right.  \label{2100}
\end{equation}%
{\it where} ${\mathbb M}$ {\it is defined in (\ref{71}), and}
\begin{equation}
\,\ell _{\varphi }(\tilde{L}_{3})=\int_{{\Omega \times Y}_{1}}h\tilde{L}%
_{3}\,\text{d}x\text{d}y=2\left\vert Y_{1}\right\vert \int_{{\omega }}h{\cal %
M}(\tilde{L}_{3})\,\text{d}\hat{x}.  \label{Eqellphi}
\end{equation}

\medskip \noindent {\bf Remark }: as usual, the convergence is a priori for
a subsequence and the uniqueness of the solution to the limit problem shows
a posteriori that the whole family converges.

\medskip \noindent {\bf Proof } First, successive multiplications of (\ref%
{63}) by $a^{2}$, $a$ and $1$ and passing to the limit yield for all ${\bf V}%
\in {\bf W}_{ad}^{1}$
\begin{equation}
\left\{
\begin{array}{l}
\dint_{{\Omega }\times Y}{\bf M}^{02}{\bf (V})\text{ }{\cal R}\text{ }{^{t}}%
{\bf M}\,\text{d}x\text{d}y=0,\medskip \  \\
\dint_{{\Omega }\times Y}{\bf M}^{01}{\bf (V})\text{ }{\cal R}\text{ }{^{t}}%
{\bf M\,}\text{d}x\text{d}y+2G\dint_{{\Omega }\times Y_{1}}{\cal M}(L_{3})%
{\cal M}\left( \partial _{3}\psi \right) \text{ d}x\text{d}y\medskip \\
\text{\ \ \ \ \ }+2G_{1}\dint_{{\Omega }\times Y_{1}}{\partial }_{\alpha }%
{\cal M}(L_{3}){\partial }_{\alpha }{\cal M}(\partial _{3}\psi )\,\text{d}x%
\text{d}y=\ell _{\varphi }(\tilde{L}_{3})\text{ if }{\bf M}^{02}{\bf (V)=M}%
^{11}{\bf (V)=0},\medskip \  \\
\dint_{{\Omega }\times Y}{\bf M}^{00}{\bf (V})\text{ }{\cal R}\text{ }{^{t}}%
{\bf M}\,\text{d}x\text{d}y=\ell _{u}({\bf v})\text{ if }{\bf M}^{02}{\bf %
(V)=M}^{11}{\bf (V)=M}^{01}{\bf (V)=M}^{10}{\bf (V)=0}.%
\end{array}%
\right.  \label{28}
\end{equation}%
\smallskip Second, let ${\bf v}=(\bar{v}_{\alpha }-x_{3}{\partial }%
_{y_{\alpha }}v_{3}^{2},v_{3}^{2},0)$ $\in $ ${\bf V}_{KL}^{1}\cap {\bf W}%
_{ad}^{1}$, with $v_{3}^{2}=0$ in a neighbourhood of ${\Gamma }_{D}$, and
choose
\[
{\bf V}^{\varepsilon }={\varepsilon }\ (\bar{v}_{1}-x_{3}{\partial }%
_{y_{_{1}}}v_{3}^{2},\bar{v}_{2}-x_{3}{\partial }_{y_{_{2}}}v_{3}^{2},{%
\varepsilon }v_{3}^{2},0)
\]%
in (\ref{63}), so that $\varepsilon {\bf M}^{00}{\bf (V}^{\varepsilon }{\bf )%
}$ tends to ${\bf 0}$, ${\bf M}^{02}{\bf (V}^{\varepsilon }{\bf )=0}$, and
\[
\frac{1}{a}{\bf M}^{01}{\bf (V}^{\varepsilon }{\bf )}+\frac{1}{a\varepsilon }%
{\bf M}^{11}{\bf (V}^{\varepsilon })=({\bf 0}_{2\times 2}{,}\dfrac{{%
\varepsilon }^{2}}{2a}{\partial }_{\alpha }v_{3}^{2},{\bf 0}_{4}).
\]%
Using lemma 7.3, passing to the limit in (\ref{63}) yields
\begin{equation}
\dint\nolimits_{{\Omega }\times Y}{\bf M}^{10}{\bf (V})\ {\cal R\ }{^{t}}%
{\bf M}\,\text{d}x\text{d}y=0\,  \label{29}
\end{equation}%
for any ${\bf V=(v,}0)\in {\bf W}_{ad}^{1}$ with ${\bf v}\in {\bf V}%
_{KL}^{1} $. Now, summing the equations (\ref{28})-(\ref{29}) and using then
the density lemma 7.1, we get (\ref{2100}).

\smallskip \noindent Last, uniqueness of the solution to this system is a
simple application of Lax-Milgram lemma$.\ \ $ \ $\Box $

\bigskip \noindent {\bf 7.2 Proof of Theorem 5.1. }

\smallskip \noindent From now on, an important use is made of the assumption
that ${\cal R}$ does not depend on $x_{3}$.

\bigskip \noindent {\bf 7.2.1 Proof of Theorem 5.1 : Dirichlet conditions. }

\smallskip \noindent We recall that here ${\cal M}(L_{3})=\varphi _{c}$ and
that ${\cal M}$\ and ${\cal N}$ have been defined in (\ref{1003}). The proof
is in two steps. The first one consists in eliminating $K_{i3}$ and ${\cal N}%
(L_{3})$; the second one is for eliminating the $y-$variable. A key for many
simplifications is that ${\cal R}$ does not depend on $x_{3}$ and so ${\cal M%
}({\bf M})$ and ${\cal N}({\bf M})$ bring separate contributions.

\smallskip \noindent {\bf First step.} Remembering that in the present case,
$G=G_{1}=0$ and $h=0$, Equation (\ref{2100}) of Lemma 7.5 reduces to%
\[
\forall {\widetilde{{\bf M}}}\in {\mathbb M},\ \ \dint\nolimits_{{\Omega }%
\times Y}{\widetilde{{\bf M}}}\ {\cal R\ }{^{t}}{\bf M}\text{ d}x\text{d}%
y=\ell _{u}({\bf v}).
\]%
Splitting ${\bf M}$ as ${\bf M}={\cal M}({\bf M})+{\cal N}({\bf M})$ and
taking advantage of the fact that ${\cal R}$ does not depend on $x_{3}$,
this is equivalent to
\begin{equation}
\forall {\widetilde{{\bf M}}}\in {\mathbb M},\ \ \dint\nolimits_{{\Omega }%
\times Y}{\cal M}({\widetilde{{\bf M}}})\ {\cal R\ }{^{t}}{\cal M}({\bf M})%
\text{ d}x\text{d}y+\dint\nolimits_{{\Omega }\times Y}{\cal N}({\widetilde{%
{\bf M}}})\ {\cal R\ }{^{t}}{\cal N}({\bf M})\,\text{d}x\text{d}y=\ell _{u}(%
{\bf v}).  \label{81}
\end{equation}%
Let us chose ${\widetilde{{\bf M}}:}={\cal M}({\widetilde{{\bf M}}})$ as
test functions in (\ref{81}) where ${\widetilde{{\bf M}}}\in {\mathbb M}%
^{-1}\oplus {\mathbb M}^{-2}$. Because the function ${\bf v}$
associated with
such a ${\widetilde{{\bf M}}}$ is ${\bf v}={\bf 0,}$ and because ${\cal %
N\circ M}=0$, we get
\[
\forall {\widetilde{{\bf M}}}\in {\mathbb M}^{-1}\oplus {\mathbb
M}^{-2},\ \
\dint\nolimits_{{\Omega }\times Y}{\cal M}({\widetilde{{\bf M}}})\ {\cal R\ }%
{^{t}}{\cal M}({\bf M})\text{ d}x\text{d}y=0.
\]%
Using again that ${\cal R}$ is independant of $x_{3}$, this may be rewritten
as
\[
\forall {\widetilde{{\bf M}}}\in {\mathbb M}^{-1}\oplus {\mathbb
M}^{-2},\ \
\dint\nolimits_{{\Omega }\times Y}{}{\widetilde{{\bf M}}\ }{\cal R\ }{^{t}}%
{\cal M}({\bf M})\text{ d}x\text{d}y=0.
\]%
This holds true in particular for ${\widetilde{{\bf M}}}$ of the form ${%
\widetilde{{\bf M}}=}({\bf 0}_{2\times 2},\tilde{K}_{{i}3},{\bf 0}_{3})$ and
thus shows that
\begin{equation}
\Pi _{2}{\cal R}\ {^{t}}{\cal M}({\bf M})={\bf 0}  \label{projMDirichlet}
\end{equation}%
On the other hand, Lemma 7.4 shows that one may decompose ${\bf M}$ as

\begin{equation}
{\bf M}=\Pi {\bf M}+{\bf M}^{0}  \label{82}
\end{equation}%
where $\Pi {\bf M}\in \phi _{c}+{\mathbb M}^{-1}\oplus {\mathbb M}^{-2}{\bf ,}$ $%
{\bf M}^{0}\in {\mathbb M}^{0}$, and that ${\cal M}(\Pi {\bf M})$
may be written as
\[
{\cal M}(\Pi {\bf M})%
\begin{tabular}[t]{l}
$=({\bf 0}_{2\times 2},{\cal M}(K_{{i}3}),{\bf 0}_{2},{\cal M}(L_{3}))$ \\
$=\phi _{c}+\Pi _{2}{\cal M}({\bf M}){\bf .}$%
\end{tabular}%
\]%
Hence, applying ${\cal M}$ to (\ref{82}) we get

\begin{equation}
{\cal M}({\bf M})%
\begin{tabular}[t]{l}
$=\phi _{c}+{\cal M}(\Pi _{2}{\bf M})+{\cal M}({\bf M}^{0})$ \\
$=\phi _{c}+\Pi _{2}{\cal M}(\Pi _{2}{\bf M})+{\cal M}({\bf M}^{0}).$%
\end{tabular}
\label{M(M)Dirichlet}
\end{equation}%
Applying now $\Pi _{2}{\cal R}$ to this identity, with (\ref{projMDirichlet}%
), we obtain then%
\[
{\bf 0}=\Pi _{2}{\cal R}\ {^{t}\phi }_{c}+(\Pi _{2}{\cal R}\Pi _{2})\ {^{t}}%
{\cal M}(\Pi _{2}{\bf M})+\Pi _{2}{\cal R}\ {^{t}}{\cal M}({\bf M}^{0}),
\]%
and therefore
\begin{equation}
{^{t}}{\cal M}(\Pi _{2}{\bf M})={\bf T}_{{\cal M}}\ {^{t}}({\phi }_{c}+{\cal %
M}({\bf M}^{0})),  \label{M(M)Dir2}
\end{equation}%
where ${\bf T}_{{\cal M}}$ has been defined in (\ref{Tenseurs}).

\smallskip \noindent Similarly, let us now consider in (\ref{81}) test
functions of the form ${\cal N}(\widetilde{{\bf M}})$ where $\widetilde{{\bf %
M}}=({\bf 0}_{4},\tilde{K}_{i3},{\bf 0}_{2},\tilde{L}_{3}),\ \tilde{K}_{i3},$
$\tilde{L}_{3}$ being given functions in $L^{2}(\Omega \times Y)$ and $L^{2}(%
{\Omega }\times Y_{1})$ respectively. As ${\cal M}({\cal N}(\tilde{L}%
_{3}))=0 $, ${\cal N}(\widetilde{{\bf M}})\in {\mathbb M}^{-1}\oplus {\mathbb M}%
^{-2}$ so that this choice is correct. Proceding as for the ${\cal M-}$%
component of ${\bf M}$, we get that%
\[
\forall {\widetilde{{\bf M}}}=({\bf 0}_{4},\tilde{K}_{i3},{\bf 0}_{2},\tilde{%
L}_{3}),\ \ \dint\nolimits_{{\Omega }\times Y}{}{\widetilde{{\bf M}}\ }{\cal %
R}\ {^{t}}{\cal N}({\bf M})\text{ d}x\text{d}y=0,
\]%
or in other \ words, $\Pi {\cal R}\ {^{t}}{\cal N}({\bf M})={\bf 0}$. Then,
from (\ref{82}),
\[
{\bf 0}=\Pi {\cal R}\ {^{t}}{\cal N}(\Pi {\bf M})+\Pi {\cal R}\ {^{t}}{\cal N%
}({\bf M}^{0}),
\]%
and thus,
\begin{equation}
{^{t}}{\cal N}(\Pi {\bf M})={\bf T}_{{\cal N}}\ {^{t}}{\cal N}({\bf M}^{0})
\label{N(M)Dir}
\end{equation}%
where ${\bf T}_{{\cal N}}$ has been defined in (\ref{Tenseurs})$.$ Thus,
with (\ref{M(M)Dirichlet}), (\ref{M(M)Dir2}), (\ref{82}) and (\ref{N(M)Dir})
we have

\begin{equation}
{\cal M}({\bf M})=(\text{Id}+{\bf T}_{{\cal M}})\ {^{t}}({\cal M}({\bf M}%
^{0})+{\phi }^{c}),\text{ and }{\cal N}({\bf M)}=(\text{Id}+{\bf T}_{{\cal N}%
})\ {^{t}}{\cal N}({\bf M}^{0}).  \label{84}
\end{equation}%
Now, if we choose - in order to maintain a symmetric form to our tensors- in
(\ref{81}) test functions of the special form ${\widetilde{{\bf M}}=}{\cal M}%
({\widetilde{{\bf M}}}^{0}){(Id}+\ {^{t}}{\bf T}_{{\cal M}})+{\cal N}({%
\widetilde{{\bf M}}}^{0}){(Id}+\ {^{t}}{\bf T}_{{\cal N}})$ with ${%
\widetilde{{\bf M}}}^{0}\in {\mathbb M}^{0}$, taking (\ref{84})\
into account, we get that ${\bf M}^{0}$ is the unique solution in
${\phi }_{c}+{\mathbb M}^{0}
$ of%
\begin{equation}
\left\{
\begin{array}{l}
\forall \,{\widetilde{{\bf M}}}\in {\mathbb M}^{0},\ \ \dint\nolimits_{{\Omega }%
\times Y}\left( {\cal M}({\widetilde{{\bf M}}})\ {\cal R}_{{\cal M}}\ {^{t}}%
{\cal M}({\bf M}^{0})+{\cal N}({\widetilde{{\bf M}}})\ {\cal R}_{{\cal N}}\ {%
^{t}}{\cal N}({\bf M}^{0})\right) \text{ d}x\text{d}y\medskip  \\
\ \ \ \ \ \ \ \ \ \ \ \ \ \ \ \ \ \ \ =\ell _{u}({\bf v})-\dint\nolimits_{{%
\Omega }\times Y}{\cal M}({\widetilde{{\bf M}}})\ {\cal R}_{{\cal M}}\ ^{t}{%
\phi }_{c}\,\text{d}x\text{d}y,%
\end{array}%
\right.   \label{85}
\end{equation}%
where ${\bf v}$ is the vector in ${\bf V}_{KL}$ associated with ${\widetilde{%
{\bf M}}}$, and ${\cal R}_{{\cal M}}$ and ${\cal R}_{{\cal N}}$\ have been
defined in (\ref{Tenseurs}).

\smallskip \noindent {\bf Second step.} It remains to eliminate ${\bf u}^{1}$%
. For that purpose, we use the usual arguments of linear homogenization.
First note that the definitions (\ref{defi KL}) of ${\bf V}_{KL}$ and (\ref%
{defiKL1}) of ${\bf V}_{KL}^{1}$ imply that
\begin{equation}
\left\{
\begin{array}{l}
\forall \,{\widetilde{{\bf M}}}\in {\mathbb M}^{0},\ \ \ {\cal M}({\widetilde{%
{\bf M}}})=\left( s_{{\alpha }\beta }({\bf \bar{v}})+S_{{\alpha }\beta }(%
{\bf \bar{v}}^{1}),{\bf 0}_{6}\right) ,\medskip \\
\forall \,{\widetilde{{\bf M}}}\in {\mathbb M}^{0},\ \ \ {\cal N}({\widetilde{%
{\bf M}}})=-x_{3}\ \left( {\partial }_{{\alpha }\beta }^{2}v_{3}+{\partial }%
_{y_{\alpha }y_{\beta }}^{2}v_{3}^{2},{\bf 0}_{6}\right) ,%
\end{array}%
\right.  \label{86}
\end{equation}%
where ${\bf v}\in {\bf V}_{KL}$ and $(\bar{v}_{1}^{1},\bar{v}%
_{2}^{1},v_{3}^{2})\in {\bf V}_{KL}^{1}$ are the fields associated with ${%
\widetilde{{\bf M}}}$.

\smallskip \noindent Also, recall that
\[
{\bf M}^{0}=((s_{{\alpha }\beta }({\bf u})+(S_{{\alpha }\beta }({\bf u}%
^{1}))_{{\alpha }\beta =1,2},{\bf 0}_{6})\text{ with }{\bf u}\in {\bf V}_{KL}%
\text{ and }{\bf u}^{1}\in {\bf V}_{KL}^{1}.
\]

\smallskip \noindent Considering ${\bf v}={\bf 0}$ in (\ref{85}), we thus
get that for all $\overline{{\bf v}}^{1}\in (H_{\sharp }^{1}(Y))^{2},$%
\begin{equation}
\dint\nolimits_{Y}S_{{\alpha }\beta }(\overline{{\bf v}}^{1})R_{{\cal M}{%
\alpha }\beta {\lambda }\mu }S_{{\lambda }\mu }({\bf \bar{u}}^{1})\text{ d}%
y=-\dint\nolimits_{Y}\,S_{{\alpha }\beta }(\overline{{\bf v}}^{1})(R_{{\cal M%
}{\alpha }\beta {\lambda }\mu }s_{{\lambda }\mu }({\bf \bar{u}})+d_{{\cal M}3%
{\alpha }\beta }\varphi _{c})\,\text{d}y,  \label{5homoDirM}
\end{equation}%
and for all $v_{3}^{2}\in H_{\sharp }^{2}(Y)$,
\begin{equation}
\dint\nolimits_{Y}{\partial }_{y_{\alpha }y_{\beta }}^{2}v_{3}^{2}\ R_{{\cal %
N}{\alpha }\beta {\lambda }\mu }{\partial }_{y_{\lambda }y_{\mu
}}^{2}u_{3}^{2}\,\text{d}y=-\dint\nolimits_{Y}\,{\partial }_{y_{\alpha
}y_{\beta }}^{2}v_{3}^{2}\ R_{{\cal N}{\alpha }\beta {\lambda }\mu }{%
\partial }_{{\lambda }\mu }^{2}u_{3}\,\text{d}y  \label{5homoDirN}
\end{equation}%
almost everywhere in $\omega $.

\smallskip \noindent Equation (\ref{5homoDirM}) implies that ${\bf \bar{u}}%
^{1}={\bf u}_{{\cal M}}^{{\rho \xi }}s_{{\rho \xi }}({\bf \bar{u}})+{\bf u}_{%
{\cal M}}^{3}\varphi _{c}$, where the {\em coefficients} functions ${\bf u}_{%
{\cal M}}^{{\rho \xi }}$, ${\bf u}_{{\cal M}}^{3}$ are defined in (\ref{24})
and (\ref{241}), and therefore%
\[
s_{{\lambda }\mu }({\bf \bar{u}})+S_{{\lambda }\mu }({\bf \bar{u}}%
^{1})=(\delta _{{\lambda }\mu ,{\rho \xi }}+S_{{\lambda }\mu }({\bf u}_{%
{\cal M}}^{\rho \xi }))s_{\rho \xi }({\bf \bar{u}})+S_{{\lambda }\mu }({\bf u%
}_{{\cal M}}^{3})\varphi _{c}\text{.}
\]

\smallskip \noindent This is for the ${\cal M}$-part of ${\bf M}^{0}$.
Similarly, Equation (\ref{5homoDirN}) implies that and $u_{3}^{2}={\bf u}_{%
{\cal N}}^{{\rho \xi }}{\partial }_{{\rho \xi }}^{2}u_{3}$ where the ${\bf u}%
_{{\cal N}}^{{\rho \xi }}$ are also defined by (\ref{24}), and this leads to

\[
{\partial }_{{\lambda }\mu }^{2}u_{3}+{\partial }_{y_{\lambda }y_{\mu
}}^{2}u_{3}^{2}=(\delta _{{\lambda }\mu ,\rho \xi }+{\partial }_{y_{\lambda
}y_{\mu }}^{2}({\bf u}_{{\cal N}}^{\rho \xi })){\partial }_{\rho \xi
}^{2}u_{3}.
\]

\smallskip \noindent Now we take in (\ref{85}) test functions as in (\ref{86}%
), with the special forms ${\bf \bar{v}}^{1}={\bf u}_{{\cal M}}^{\gamma
\delta }s_{{\gamma }{\delta }}(\overline{{\bf v}})$ and $v_{3}^{2}={\bf u}_{%
{\cal N}}^{\gamma \delta }{\partial }_{{\gamma }{\delta }}^{2}v_{3}$. This
leads to the system (\ref{EffectifDirichlet}) with
\[
\left\{
\begin{array}{l}
d_{{\cal M}3\gamma \delta }^{H}=\dint\nolimits_{Y}(\delta _{{\alpha }\beta
,\gamma \delta }+S_{{\alpha \beta }}({\bf u}_{{\cal M}}^{{\gamma }{\delta }%
}))\left( R_{{\cal M}{\alpha }\beta {\lambda }\mu }S_{{\lambda }\mu }({\bf u}%
_{{\cal M}}^{3})+d_{{\cal M}3\alpha \beta }\right) \text{ d}y\text{,}\medskip
\\
R_{{\cal M}\gamma \delta \rho \xi }^{H}=\dint\nolimits_{Y}(\delta _{{\alpha }%
\beta ,\gamma \delta }+S_{{\alpha \beta }}({\bf u}_{{\cal M}}^{{\gamma }{%
\delta }}))R_{{\cal M}{\alpha }\beta {\lambda }\mu }(\delta _{{\lambda }\mu
,\rho \xi }+S_{{\lambda }\mu }({\bf u}_{{\cal M}}^{\rho \xi }))\text{ d}%
y,\medskip \\
R_{{\cal N}\gamma \delta \rho \xi }^{H}=\dint\nolimits_{Y}(\delta _{{\alpha }%
\beta ,\gamma \delta }+{\partial }_{y_{\alpha }y_{\beta }}^{2}u_{{\cal N}}^{{%
\gamma }{\delta }})R_{{\cal N}{\alpha }\beta {\lambda }\mu }(\delta _{{%
\lambda }\mu ,\rho \xi }+{\partial }_{y_{{\lambda }}y_{\mu }}^{2}u_{{\cal N}%
}^{\rho \xi })\text{ d}y.%
\end{array}%
\right.
\]%
which are exactly the tensors announced in formulae (\ref{TenseurRN})-(\ref%
{TenseursM}). Note in particular that the coefficient $2$ appears on the
first term on left hand side of (\ref{EffectifDirichlet}) because $%
\dint\nolimits_{-1}^{1}$d$x_{3}=2$ while $\dfrac{2}{3}=\dint%
\nolimits_{-1}^{1}x_{3}^{\ 2}$d$x_{3}$ appears on the second term. \ \ \ \ \
{$\Box $}\bigskip

\bigskip \noindent {\bf 7.2.2 Proof of Theorem 5.1 : nonlocal mixed
conditions. \ }

\medskip \noindent The proof is very closed to the one for Dirichlet
conditions. The main difference is that in the Dirichlet case, ${\cal M}%
(L_{3})=\varphi _{c}$ is a given data, while here ${\cal M}(L_{3})$ is an
unknown function which can not be eliminated.

\medskip \noindent Let us eliminate the other transverse components $%
K_{\alpha 3},K_{33}$ and ${\cal N}(L_{3}).$

\medskip \noindent The weak formulation (\ref{2100}) of Lemma 7.5, implies
now that
\begin{equation}
\left\{
\begin{array}{l}
\forall {\widetilde{{\bf M}}}\in {\mathbb M},\ \ \dint\nolimits_{{\Omega }%
\times Y}\left( {\cal M}({\widetilde{{\bf M}}})\ {\cal R\ }^{t}{\cal M}({\bf %
M})\,+{\cal N}({\widetilde{{\bf M}}})\ {\cal R\ }^{t}{\cal N}({\bf M}%
)\right) \,\text{d}x\text{d}y\medskip  \\
\ \ \ \ +2\dint\nolimits_{{\Omega }\times Y_{1}}\left( G{\cal M}(L_{3}){\cal %
M}(\tilde{L}_{3})+G_{1}{\partial }_{\alpha }{\cal M}(L_{3}){\partial }%
_{\alpha }{\cal M}(\tilde{L}_{3})\right) \text{ d}x\text{d}y=\ell _{u}({\bf v%
})+\ell _{\varphi }(\tilde{L}_{3}).%
\end{array}%
\right.   \label{8100}
\end{equation}

\medskip \noindent For ${\cal M}({\bf M})$, we use test functions such that $%
\tilde{L}_{3}=0,$ as in Section 7.2.1, and thus the computation go on the
same way and the result is the same, except that $\varphi _{c}$ is now
replaced by ${\cal M}(L_{3})$. For ${\cal N}({\bf M})$, as ${\cal M\circ N}$
and as $\ell _{\varphi }({\cal N}(\tilde{L}_{3}))=0$ (see (\ref{Eqellphi}))
the terms in $G,$ $G_{1}$ and $\ell _{\varphi }$ do not play any part as
well, and the computation also go as in 7.2.1.\ We thus get

\[
{\cal M}({\bf M})=(\text{Id}+{\bf T}_{{\cal M}})\ {^{t}}\left( {\cal M}({\bf %
M}^{0})+{\Lambda }_{3}\right) ,\text{ and }{\cal N}({\bf M)}=(\text{Id}+{\bf %
T}_{{\cal N}})\ {^{t}}{\cal N}({\bf M}^{0}).
\]%
where ${\Lambda }_{3}=({\bf 0}_{9},{\cal M}(L_{3}))$. Then, with a suitable
choice of test functions, (\ref{8100}) implies that $({\bf M}^{0},{\cal M}({L%
}_{3}))$ is the unique solution in\ ${\mathbb M}^{0}\times {L}^{2}(\omega )$ of%
\begin{equation}
\left\{
\begin{array}{l}
\forall \,\widetilde{{\bf M}}^{0}\in {\mathbb M}^{0}\text{, }\forall \tilde{L}%
_{3}{\in L}^{2}(\omega )\text{,} \\
\ \ \ \ \dint\nolimits_{{\Omega }\times Y}\left( ({\cal M}(\widetilde{{\bf M}%
}^{0})+\widetilde{\Lambda }_{3})\ {\cal R}_{{\cal M}}\ ^{t}\left( {\cal M}(%
{\bf M}^{0})+{\Lambda }_{3}\right) +{\cal N}({\widetilde{{\bf M}}}^{0})\
{\cal R}_{{\cal N}}\ ^{t}{\cal N}({\bf M}^{0})\right) \,\text{d}x\text{d}%
y\medskip \\
\ \ \ \ \ \ +2\dint\nolimits_{{\Omega }\times Y_{1}}\left( G{\cal M}(L_{3})%
\tilde{L}_{3}\,+G_{1}{\partial }_{\alpha }{\cal M}(L_{3}){\partial }_{\alpha
}\tilde{L}_{3}\right) \,\text{d}x\text{d}y=\ell _{u}({\bf v})+\ell _{\varphi
}(\tilde{L}_{3})\,%
\end{array}%
\right.  \label{8500}
\end{equation}%
where $\widetilde{\Lambda }_{3}=({\bf 0}_{9},\tilde{L}_{3})$ and ${\bf v}$
is the vector of ${\bf V}_{KL}$ associated with $\,\widetilde{{\bf M}}^{0}$.

\medskip \noindent To eliminate the local variable $y$, we proceed as in the
Dirichlet case, with $\widetilde{\Lambda }_{3}=({\bf 0}_{9},\tilde{L}_{3})$
instead of ${\phi }_{c}$. First, letting $L_{3}^{0}={\cal M(}L_{3}{\cal )}$
and considering in (\ref{8500}) test functions such that $\tilde{L}_{3}=0$
and ${\bf \bar{v}=0}$, we get that ${\bf \bar{u}}^{1}={\bf u}_{{\cal M}}^{{%
\rho \xi }}s_{{\rho \xi }}({\bf \bar{u}})+{\bf u}_{{\cal M}}^{3}L_{3}^{0}$, $%
u_{3}^{2}={\bf u}_{{\cal N}}^{{\rho \xi }}{\partial }_{{\rho \xi }}^{2}u_{3}$%
. Then, remarking that fot all$\,\widetilde{{\bf M}}^{0}\in {\mathbb
M}^{0}$
and all $\tilde{L}_{3}{\in L}^{2}(\omega ),$%
\[
\left\{ \text{%
\begin{tabular}{l}
$({\cal M}(\widetilde{{\bf M}}^{0})+\widetilde{\Lambda }_{3})=((s_{{\alpha }%
\beta }({\bf \bar{v}})+S_{{\alpha }\beta }({\bf \bar{v}}^{1}))_{{\alpha ,}%
\beta =1,2},{\bf 0}_{5},\tilde{L}_{3}),\medskip $ \\
${\cal N}({\widetilde{{\bf M}}}^{0})=-x_{3}\text{ }(({\partial }_{{\alpha }%
\beta }^{2}v_{3}+{\partial }_{y_{\alpha }y_{\beta }}^{2}v_{3}^{2})_{{\alpha ,%
}\beta =1,2},{\bf 0}_{6}),$%
\end{tabular}%
}\right.
\]%
with a suitable choice of test functions in (\ref{8500}) we obtain the
announced model. $\ \ \ \Box $

\bigskip \noindent Note that the above calculations would be more
complicated with non metallized inclusions, as ${\cal M}(L_{3})$ would then
depend on $y$.

\bigskip \noindent {\bf 7.2.3 Proof of Theorem 5.1 : local mixed conditions.
}

\smallskip \noindent Here $G_{1}=0$ so that the effective equations
simplifies into%
\[
\left\{
\begin{array}{l}
\dint_{{\omega }}\left( 2(s_{{\alpha }\beta }({\bf \bar{v}}),\tilde{L}%
_{3})\left(
\begin{array}{cc}
R_{{\cal M}{\alpha }\beta {\gamma }{\delta }}^{H} & d_{{\cal M}3{\alpha }%
\beta }^{H} \\
e_{{\cal M}3{\gamma }{\delta }}^{H} & c_{{\cal M}33}^{H}+2|Y_{1}|G%
\end{array}%
\right) \left(
\begin{array}{c}
s_{{\gamma \delta }}({\bf \bar{u}}) \\
L_{3}^{0}%
\end{array}%
\right) \right) \text{ d}\hat{x}\medskip  \\
\ \ \ \ \ \ \ \ \ \ \ \ \ \ \ \ \ \ \ \ \ \ \ \ \ \ \ \ \ +\dfrac{2}{3}%
\dint_{{\omega }}{\partial }_{{\alpha }\beta }^{2}v_{3}R_{{\cal N}{\alpha }%
\beta {\gamma \delta }}^{H}{\partial }_{{\gamma \delta }}^{2}u_{3}\text{ d}%
\hat{x}=\ell _{u}({\bf v})+2\left\vert Y_{1}\right\vert \dint_{{\omega }}%
\tilde{L}_{3}\text{ }h\text{ d}\hat{x}.\medskip
\end{array}%
\right.
\]%
This simplification allows to eliminate the unknown $L_{3}^{0}$ that can be
computed explicitly in terms of $s_{{\gamma \delta }}({\bf \bar{u}})$.\ We
choose ${\bf v=0}$ then
\[
\forall \tilde{L}_{3}\in L^{2}(\omega )\text{, }\dint\nolimits_{{\omega }%
}\left( e_{{\cal M}3{\gamma }{\delta }}^{H}s_{{\gamma \delta }}({\bf \bar{u}}%
)+\left( c_{{\cal M}33}^{H}+2|Y_{1}|G\right) L_{3}^{0}-\left\vert
Y_{1}\right\vert h\right) \tilde{L}_{3}\text{ d}\hat{x}=0
\]%
so that
\[
\left( c_{{\cal M}33}^{H}+2|Y_{1}|G\right) L_{3}^{0}=\left\vert
Y_{1}\right\vert h-e_{{\cal M}3{\gamma }{\delta }}^{H}s_{{\gamma \delta }}(%
{\bf \bar{u}})\text{ a.e. in }\omega \text{.}
\]%
Now replacing $L_{3}^{0}$ and \ restricting ourselves to test functions with
$\tilde{L}_{3}=0$, we get the announced model, and thus conclude the proof.
\ \ $\Box $

\bigskip \noindent {\bf 8. Proof of Theorem 6.1. }

\smallskip \noindent This section is devoted to the derivation of Theorem
6.1. The proof is based on the general results of Lemma 4.1. As
homogenization and plate theory act here on the same level, the proof is in
two steps only: first, characterization of the limit ${\bf M}$ defined in
Lemma 4.1; this is the aim of Section 8.1; second, elimination of the local
variable $(y_{1},y_{2},x_{3})$; this is achieved in Section 8.2.\bigskip

\bigskip \noindent {\bf 8.1 Step 1: Characterization of the limit }${\bf M}$%
{\bf . }

\medskip \noindent {\bf 8.1.1 Characterization of }${\bf M}$.

\medskip \noindent {\bf Lemma 8.1.{\rm {\it \ }}}{\it Let }${\bf M}={\bf (K,L%
})${\it \ be the limit, up to a subsequence, of }$\left( {\bf M}^{a}({\bf U}%
^{b})\right) $, t{\it hen}\newline
{\it (i)}\ {\it there exist }${\bf u\in V}_{KL}$ {\it and} ${\bf \hat{u}}%
^{1}\in (L^{2}(\Omega ;H_{\sharp }^{1}(Y))\cap L^{2}(\omega ;H^{1}(Z))^{3}$
{\it such that}%
\[
\begin{tabular}{l}
$\forall \alpha ,\beta \in \left\{ 1,2\right\} $,$\ K_{a\beta }=s_{\alpha
\beta }({\bf u})+S_{\alpha \beta }({\bf \hat{u}}^{1}),\medskip $ \\
$\forall i\in \left\{ 1,2,3\right\} $,$\ K_{i3}=S_{i3}({\bf \hat{u}}^{1});$%
\end{tabular}%
\]%
{\it (ii) there exists }$\varphi ^{1}$ $\in $ $L^{2}({\omega };H^{1}(Z_{1}))$%
{\it \ such that }$L=\nabla _{z}\varphi ^{1}$;$\medskip $\newline
{\it (iii) in the case of Dirichlet conditions, one may chose} $\varphi ^{1}$
{\it so that} $\varphi ^{1}=x_{3}\varphi _{c}+\hat{\varphi}^{1}$ {\it with} $%
\hat{\varphi}^{1}\in L^{2}({\omega };H^{1}(Z_{1}))$ {\it and }$\hat{\varphi}%
^{1}=0$ {\it on} $\Gamma ^{+}\cup \Gamma ^{-}$;$\medskip $\newline
{\it (iv) in the case of mixed conditions one may chose }$\varphi ^{1}$ {\it %
so that} $\varphi ^{1}=\left( 1+x_{3}\right) {\cal M}(L_{3})+\hat{\varphi}%
^{1}$ {\it with }$\hat{\varphi}^{1}\in L^{2}({\omega };H^{1}(Z_{1}))$ {\it %
and} $\hat{\varphi}^{1}=0$ {\it on} $\Gamma ^{+}\cup \Gamma ^{-}$.

\medskip \noindent {\bf Proof :{\rm {\it \ }}}First, from Lemma 4.1 we know
that there is $({\bf u,u}^{1})\in V_{KL}\times L^{2}(\Omega ;{\bf H}_{\sharp
}^{1}(Y)/{\mathbb R})^{2}$ such that $K_{\alpha \beta }=s_{{\alpha }\beta }{\bf %
(u})+S_{{\alpha }\beta }{\bf (u}^{1})$. Then here, passing to the limit in (%
\ref{eq u un}) as $a/\varepsilon \rightarrow 1$, we get, taking ${\bf u}\in
V_{KL}$ into account, that for all $v\in {\cal D}({\Omega };{\cal C}_{\sharp
}^{\infty }(Y)),$%
\begin{equation}
\dint\nolimits_{{\Omega }\times Y}({\partial }_{y_{1}}u_{2}^{1}-{\partial }%
_{y_{2}}u_{1}^{1}){\partial }_{3}v\text{ d}x\text{d}y=2\dint\nolimits_{{%
\Omega }\times Y}(K_{23}{\partial }_{y_{1}}v-K_{{1}3}{\partial }_{y_{2}}v)%
\text{ d}x\text{d}y\text{.}  \label{35}
\end{equation}%
Hence, see for instance [24, Theorem 5], there is $u_{3}^{2}$ and $c_{1},$ $%
c_{2}$ not depending on $y$ such that for ${\alpha }=1,2$, ${\partial }%
_{3}u_{\alpha }^{1}-{2}K_{{\alpha }3}=c_{\alpha }-{\partial }_{y_{\alpha
}}u_{3}^{2}$. Also, note that $u_{3}^{2}$ is defined up to the addition of a
function of $x.$ Here, as ${\bf u}^{1}$ is also defined up to the addition
of a function of $x$, one may choose $u_{\alpha }^{1}$ so that $c_{\alpha
}=0 $. Letting ${\bf \hat{u}}^{1}=\left(
u_{1}^{1},u_{2}^{2},u_{3}^{2}\right) ,$ then $K_{\alpha 3}=S_{{\alpha }3}(%
\widehat{{\bf u}}^{1})$.

\medskip \noindent We now prove that $K_{33}={\partial }_{3}u_{3}^{2}.$ With
a few integrations by parts, one easily gets that for any $v\in {\cal D}({%
\Omega }\times Y)$,%
\[
\begin{tabular}{l}
$\dint\nolimits_{\Omega }{\partial }_{\alpha }u_{\beta }^{b}{\partial }%
_{33}^{2}v^{\varepsilon }\text{d}x+\dint\nolimits_{\Omega }{\partial }%
_{\beta }u_{3}^{b}\ \left( {\partial }_{\alpha 3}^{2}v^{\varepsilon }+{%
\dfrac{1}{{\varepsilon }}}{\partial }_{3y_{\alpha }}^{2}v^{\varepsilon
}\right) \,\text{d}x=$ \\
$\,\,\,\,\,\,\,\,\,\,\,\,\,\,\,\,\,\,\,\,\,\,\,\,\,\,\,\,\,\,\,\,\,2a\dint%
\nolimits_{\Omega }K_{\beta 3}^{a}({\bf u}^{b}){\partial _{\alpha 3}^{2}}%
v^{\varepsilon }\,\,\text{d}x+{2}\dint\nolimits_{\Omega }K_{\beta 3}^{a}(%
{\bf u}^{b}){\partial }_{3y_{\alpha }}^{2}v^{\varepsilon }\,\text{d}x.$%
\end{tabular}%
\]%
Besides, from Lemma 4.1, we know that $\left( a^{2}K_{33}^{a}({\bf u}%
^{b})\right) $ and $\left( aK_{33}^{a}({\bf u}^{b})\right) $ tend to zero
and that $\left( K_{33}^{a}({\bf u}^{b})\right) $ tends to $K_{33}$. Thus,
integrating by parts again (for the second term on the left hand side) and
passing to the limit, we get%
\[
\dint\nolimits_{{\Omega }\times Y}({\partial }_{\alpha }u_{\beta }+{\partial
}_{y_{\alpha }}u_{\beta }^{1}){\partial }_{33}^{2}v\,\text{d}x\text{d}%
y+\dint\nolimits_{{\Omega }\times Y}K_{33}{\partial }_{y_{\alpha }y_{\beta
}}^{2}v\,\text{d}x\text{d}y={2}\dint\nolimits_{{\Omega }\times Y}K_{\beta 3}{%
\partial }_{y_{\alpha }3}^{2}v\,\text{d}x\text{d}y
\]%
for all $v\in {\cal D}({\Omega }\times Y)$. But because ${\bf u}\in V_{KL}$,
this is equivalent to
\begin{equation}
\dint\nolimits_{{\Omega }\times Y}{\partial }_{y_{\alpha }}u_{\beta }^{1}{%
\partial }_{33}^{2}v\,\text{d}x\text{d}y+\dint\nolimits_{{\Omega }\times
Y}K_{33}{\partial }_{y_{\alpha }y_{\beta }}^{2}v\,\text{d}x\text{d}y={2}%
\dint\nolimits_{{\Omega }\times Y}K_{\beta 3}{\partial }_{3y_{\alpha
}}^{2}v\,\text{d}x\text{d}y.  \label{36}
\end{equation}%
\smallskip \noindent Choosing $v=w_{\alpha }$ in (\ref{36}), summing for ${%
\alpha }=1$ to $2$ and using the fact that for each $v$ in ${\cal D}({\Omega
}\times Y)$, there is a ${\bf w}$ in $({\cal D}({\Omega }\times Y))^{2}$
such that $div_{y}\,{\bf w}=v$, (\ref{36}) implies that ${\partial }%
_{33}^{2}u_{\beta }^{1}+{\partial }_{y_{\beta }}K_{33}$ $-{2\partial }%
_{3}K_{\beta 3}=0$ for $\beta =1,2$. As ${2}K_{\beta 3}={\partial }%
_{3}u_{\beta }^{1}+{\partial }_{y_{\beta }}u_{3}^{2}$, this is in turn
equivalent to ${\partial }_{y_{\beta }}({\partial }_{3}u_{3}^{2}-K_{33})=0$.
But because $u_{3}^{2}$ is defined up to an additive function of $x$ we may
choose $u_{3}^{2}$ such that $K_{33}={\partial }_{3}u_{3}^{2}.$ Then $S_{ij}(%
{\bf \hat{u}}^{1})$ $\in L^{2}({\Omega }\times Y)$ for each pair $(i,j)\in
\{1..3\}^{2}$ and therefore ${\bf \hat{u}}^{1}\in {\bf L}^{2}({\omega }%
;H^{1}(Z)).$ Last, using (\ref{35}), $u_{3}^{2}$ is $Y$-periodic. This ends
the proof of point (i).

\smallskip \noindent Now, we prove (ii) and (iii) in the case of Dirichlet
conditions. From Lemma 4.1, we already know that there exists $\varphi ^{1}$
such that $L_{1}=\partial _{y_{1}}\varphi ^{1}$ and $L_{2}=\partial
_{y_{2}}\varphi ^{1}$. Besides, passing to the limit in (\ref{Eq Lim L3})
yields%
\[
\int\nolimits_{{\Omega }\times Y_{1}}L_{3}{\partial }_{y_{\alpha }}\psi \,%
\text{d}x\text{d}y=\int\nolimits_{{\Omega }\times Y_{1}}{\partial }%
_{y_{\alpha }}\varphi ^{1}\,{\partial }_{3}{\psi }\,\text{d}x\ \text{d}y
\]%
\smallskip \noindent for all ${\alpha }=1,2$ and all ${\psi }\in {\cal D}({%
\Omega }\times Y_{1})$. Thus, ${\partial }_{3}\varphi ^{1}-L_{3}$ does not
depend on $y$ in $Y_{1}$. Hence, since $\varphi ^{1}$ is defined up to a
function of $x$, one may choose $\varphi ^{1}$ so that $L_{3}={\partial }%
_{3}\varphi ^{1}$\thinspace and we have in addition that $\varphi ^{1}\in
L^{2}({\omega };H^{1}(Z_{1}))$, so (ii) is proven. On the other hand, as we
know from Lemma 4.1 that ${\cal M}(L_{3})=\varphi _{c}$, we have that $%
\varphi _{/\Gamma ^{+}}^{1}-\varphi _{/\Gamma ^{-}}^{1}=2\varphi _{c}$.
Since so far, $\varphi ^{1}$ is still defined up to a function of $\hat{x}$,
we may choose $\varphi ^{1}$ so that $\varphi ^{1}=\varphi _{c}$ on $\Gamma
^{+}$ and $\varphi ^{1}=-\varphi _{c}$ on $\Gamma ^{-}$, or in other words,
choose $\varphi ^{1}$ as stated in (iii).

\medskip \noindent The case of mixed conditions is more complicated. Let us
reconsider the limit of $\left( {\bf L}^{a}(\varphi ^{b})\right) $ globally.
We remark that $\partial _{3}\varphi ^{b}=\partial _{3}\left( \varphi
^{b}-\varphi _{m}^{b}\right) $, and also, due to the assumption of
metallization, that $\partial _{\alpha }\varphi ^{b}=\partial _{\alpha
}\left( \varphi ^{b}-\varphi _{m}^{b}\right) $ for $\alpha =1,2$. The limit $%
{\bf L}$ of $\left( {\bf L}^{a}(\varphi ^{b})\right) $ is therefore also the
limit of $\left( {\bf L}^{a}(\varphi ^{b}-\varphi _{m}^{b})\right) $.
Besides, as
\[
\varepsilon ^{-1}\partial _{3}\left( \varphi ^{b}-\varphi _{m}^{b}\right) =%
\dfrac{a}{\varepsilon }L_{3}^{a}\left( \varphi ^{b}-\varphi _{m}^{b}\right) =%
\dfrac{a}{\varepsilon }L_{3}^{a}\left( \varphi ^{b}\right) ,
\]%
and as $\varepsilon /a\rightarrow 1$, $\left( \varepsilon ^{-1}\partial
_{3}\left( \varphi ^{b}-\varphi _{m}^{b}\right) \right) $ is bounded in $%
L^{2}\left( \Omega _{1}^{\varepsilon }\right) $ and also two-scale converge
to $L_{3}$. Thus, if we pose
\[
{\bf \tilde{L}}^{b}=\left( \partial _{\alpha }(\varphi ^{b}-\varphi
_{m}),\varepsilon ^{-1}\partial _{3}(\varphi ^{b}-\varphi _{m})\right) ,
\]%
then $\left( {\bf \tilde{L}}^{b}\right) $ two-scale converges to ${\bf L}$.

\medskip \noindent Also, by Poincar\'{e} inequality, the fact that $\left(
\varepsilon ^{-1}\partial _{3}\left( \varphi ^{b}-\varphi _{m}^{b}\right)
\right) $ is bounded in $L^{2}\left( \Omega _{1}^{\varepsilon }\right) $ and
that $\varphi ^{b}-\varphi _{m}^{b}$ vanishes on $\Gamma _{1}^{-}$implies
that $\left( \varepsilon ^{-1}\left( \varphi ^{b}-\varphi _{m}^{b}\right)
\right) $ is bounded in $L^{2}\left( \omega _{1}^{\varepsilon
};H^{1}(-1,1)\right) .$ Hence, there exists $\varphi ^{1}\in L^{2}(\omega
\times Y_{1};H^{1}(-1,1))$ such that $\left( \varepsilon ^{-1}\left( \varphi
^{b}-\varphi _{m}^{b}\right) \right) $ two-scale converges to $\varphi ^{1}$%
, and $\left( \varepsilon ^{-1}\partial _{3}\left( \varphi ^{b}-\varphi
_{m}^{b}\right) \right) $ two-scale converges to $\partial _{3}\varphi ^{1}$%
. Then, passing to the limit in%
\begin{eqnarray*}
\int\nolimits_{\Omega _{1}^{\varepsilon }}{\bf \tilde{L}}^{b}{\bf \psi }%
^{\varepsilon }\text{d}x &=&-\int\nolimits_{\Omega _{1}^{\varepsilon
}}\left( \varphi ^{b}-\varphi _{m}\right) \left( \partial _{1}\psi
_{1}^{\varepsilon }+\partial _{2}\psi _{2}^{\varepsilon }\right) \text{ d}x
\\
&&-\dfrac{1}{\varepsilon }\int\nolimits_{\Omega _{1}^{\varepsilon }}\left(
\varphi ^{b}-\varphi _{m}\right) \limfunc{div}{}_{z}{\bf \psi }^{\varepsilon
}\text{d}x+\int_{\Gamma _{1}^{\varepsilon }}\left( \varphi ^{b}-\varphi
_{m}\right) {\bf \psi }^{\varepsilon }.{\bf n}\text{ d}\sigma (x),
\end{eqnarray*}%
we get that for any function ${\bf \psi }$ in ${\em D}(\Omega \times Y_{1})$%
,
\[
\int\nolimits_{\Omega \times Y_{1}}{\bf L\psi }\text{d}x\text{d}%
y=-\int\nolimits_{\Omega \times Y_{1}}\varphi ^{1}\limfunc{div}{}_{z}{\bf %
\psi }\text{ d}x\text{d}y\text{.}
\]%
This proves that $\varphi ^{1}$ actually belongs to $L^{2}(\omega
;H^{1}(Y_{1}\times ]-1,1[))$ and that ${\bf L=\nabla }_{z}\varphi ^{1}$.

\medskip \noindent Moreover, the continuity of the trace function implies
that $\varphi ^{1}=0$ on $\Gamma ^{-}$ and, with the assumption of
metallization, that $\varphi ^{1}$\ does not depend on $y$ on $\Gamma ^{+}$%
.\ Then passing to the limit in
\[
\forall \psi \in D(\omega \times Y_{1}),\ \dfrac{a}{\varepsilon }%
\int\nolimits_{\Omega _{1}^{\varepsilon }}L_{3}^{a}\left( \varphi
^{b}\right) \psi ^{\varepsilon }\text{d}x=\dfrac{1}{2}\int\nolimits_{\omega
_{1}^{\varepsilon }}\text{tr}_{\Gamma ^{+}}\left( \frac{\varphi ^{b}-\varphi
_{m}^{b}}{\varepsilon }\right) \psi ^{\varepsilon }\text{d}\hat{x}
\]%
we get that $\varphi ^{1}=2{\cal M(}L_{3}{\cal )}$ on $\Gamma ^{+}$.\ Thus,
letting $\hat{\varphi}^{1}=\varphi ^{1}-(1+x_{3}){\cal M(}L_{3}{\cal )}$, we
get the announced result. \ \ \ $\Box $\bigskip

\medskip \noindent {\bf 8.1.2 Intermediate limit model}. We state now a
first limit model, including both global and local variables, that is
directly deduced from the convergence results of Lemma 4.1 and the
characterization in Lemma 8.1.

\medskip \noindent The limit ${\bf M}$ of $\left( {\bf M}^{a}({\bf U}%
^{b})\right) $ is completely determined by the {\em macroscopic fields }$%
{\bf U=}({\bf u,}{\cal M}(L_{3}))$ $\in ({\bf 0}_{3},\varphi _{c})+{\bf W}%
^{0}$, where we have let
\[
\begin{tabular}{l}
${\bf W}^{0}={\bf V}_{KL}\times \left\{ 0\right\} \ $in the case of
Dirichlet conditions, \\
${\bf W}^{0}={\bf V}_{KL}\times L^{2}(\omega )$ otherwise,%
\end{tabular}%
\]%
and the {\em microscopic fields }${\bf U}^{1}=\left( {\bf \hat{u}}^{1},\hat{%
\varphi}^{1}\right) \in {\bf L}^{2}(\omega ;{\bf W}^{1})$ where ${\bf W}^{1}$
is defined in (\ref{def espaces locaux}). More specifically, since $\partial
_{y_{\alpha }}\hat{\varphi}^{1}=\partial _{y_{\alpha }}\varphi ^{1}$, the
limit ${\bf M}$ takes the form ${\bf M=M}^{0}({\bf U})+{\bf M}^{1}({\bf U}%
^{1})$ where
\[
{\bf M}^{0}({\bf U})={\bf ((}s_{{\alpha }\beta }{\bf (u}))_{{\alpha },\beta
=1,2},{\bf 0}_{5},{\cal M(}L_{3}{\cal )}),
\]%
and the operator ${\bf M}^{1}$ has been defined in (\ref{Defi M1 et Sz}).

\medskip \noindent Now, we state the following intermediate limit model
using the definitions (\ref{ell u}), (\ref{Eqellphi}) of $\ell _{u}$, $\ell
_{\varphi }$ and $L_{3}^{0}={\cal M}(L_{3})$.

\medskip \noindent {\bf Lemma 8.2.{\rm {\it \ }}}{\it The limit} ${\bf M}$
{\it of} ${\bf M}^{a}({\bf U}^{b})${\bf {\rm {\it \ }}}{\it takes the form }$%
{\bf M=M}^{0}{\bf (U)+M}^{1}{\bf (U}^{1})$ {\it where }${\bf U=}({\bf u,}%
L_{3}^{0})\in ({\bf 0}_{3},\varphi _{c})+{\bf W}^{0}$ {\it and} ${\bf U}%
^{1}=\left( {\bf \hat{u}}^{1},\hat{\varphi}^{1}\right) \in {\bf L}%
^{2}(\omega ;{\bf W}^{1})$ {\it is the unique solution of this form of}%
\begin{equation}
\left\{
\begin{array}{l}
\forall {\bf V=(v,}\widetilde{L}_{3})\in {\bf W}^{0},\,\,\forall {\bf V}^{1}%
{\bf =(}\widehat{{\bf v}}^{1},{\psi }^{1})\in {\bf L}^{2}(\omega ;{\bf W}%
^{1}),\medskip \\
\dint\nolimits_{{\Omega }\times Y}({\bf M}^{0}{\bf (V)+M}^{1}{\bf (V}^{1}))%
{\cal R\ }{^{t}}({\bf M}^{0}{\bf (U)+M}^{1}{\bf (U}^{1}))\,\,\text{d}x\text{d%
}y+\medskip \\
4\left\vert Y_{1}\right\vert \dint_{{\omega }}(G\widetilde{L}%
_{3}L_{3}^{0}\,+G_{1}{\partial }_{\alpha }\widetilde{L}_{3}{\partial }%
_{\alpha }L_{3}^{0})\,\text{d}\hat{x}=l_{u}({\bf v})+\ell _{\varphi }(%
\widetilde{L}_{3}).%
\end{array}%
\right.  \label{93}
\end{equation}

\medskip \noindent {\bf Proof:{\rm {\it \ }}}Using the definitions (\ref%
{W1ad}) of ${\bf W}_{ad}^{1}$ and the notations (\ref{1331})-(\ref{2331}),
multiplying (\ref{63}) successively by ${\varepsilon }^{2}$, ${\varepsilon }$
and 1, and passing to the limit, one gets%
\[
\left\{
\begin{array}{l}
\forall {\bf V}^{1}\in {\bf W}_{ad}^{1},\ \dint\nolimits_{{\Omega }\times Y}(%
{\bf M}^{11}({\bf V}^{1})+{\bf M}^{02}({\bf V}^{1})){\cal R}\,^{t}{\bf M}\,%
\text{d}x\text{d}y=0,\medskip \\
\forall {\bf V}^{1}\in {\bf W}_{ad}^{1}\,\ \text{with }\,v_{3}^{1}=0,\
\dint\nolimits_{{\Omega }\times Y}({\bf M}^{01}({\bf V}^{1})+{\bf M}^{10}(%
{\bf V}^{1})){\cal R}\,{^{t}}{\bf M}\,\text{d}x\text{d}y\medskip \\
\ \ \ \ \ \ \ \ \ \ \ \ +2\dint\nolimits_{{\Omega }\times {Y}_{1}}(G{\cal M(}%
\widetilde{L}_{3}{\cal )}L_{3}^{0}\,+G_{1}{\partial }_{\alpha }{\cal M(}%
\widetilde{L}_{3}{\cal )}{\partial }_{\alpha }L_{3}^{0})\,\text{d}x\text{d}%
y=\ell _{\varphi }(\widetilde{L}_{3})\medskip \\
\forall {\bf v\in V}_{KL}\text{, \ }\int_{{\Omega }\times Y}{^{t}}{\bf M}%
^{00}({\bf (v,}0)){\cal R}\,{^{t}}{\bf M\,}\text{d}x\text{d}y=\ell _{u}({\bf %
v}).%
\end{array}%
\right.
\]%
Unlike the preceding models, the first two equations are coupled.
Consequently, the computation of $K_{33}$ on one hand and of $(K_{{1}3},K_{{2%
}3})$ on the other hand cannot be carried out independently. Summing up the
three equations above and using usual density results, and restricting
ourselves to test functions such that $\widetilde{L}_{3}$ does not depend on
$x_{3}$, we get the announced weak formulation. Uniqueness of the solution
follows from Lax-Milgram Lemma, and the convergence of the whole family from
the uniqueness of the solution. \ $\Box $

\bigskip \noindent {\bf 8.2 Step 2: Derivation of the models. }In this
section, we complete the proof of Theorem 6.1 by eliminating the local
variable $z$ and the corresponding unknown ${\bf U}^{1}$. We use notations
already introduced in Section 6.1.

\medskip \noindent For ${\bf W}^{1}$ defined in (\ref{def espaces locaux}),
considering (\ref{93}) with ${\bf V=0}$, we easily get%
\begin{equation}
\left\{
\begin{array}{l}
\forall {\bf V}^{1}{\bf =(v}^{1},{\psi }^{1})\in {\bf W}^{1},\medskip \\
\dint\nolimits_{Z}{\bf M}^{1}{\bf (V}^{1}){\cal R}\text{ }{^{t}}{\bf M}^{1}%
{\bf (U}^{1})\text{ d}z=-\dint\nolimits_{Z}{\bf M}^{1}({\bf V}^{1}){\cal R}%
\text{ }^{t}{\bf M}^{0}{\bf (U)}\text{d}z{\bf .}%
\end{array}%
\right.  \label{521}
\end{equation}

\medskip \noindent Let then ${\bf M}_{{\cal M}}^{0}{\bf (U)=}\left(
(s_{\gamma \delta }({\bf \bar{u})})_{\gamma ,\delta =1,2},{\bf 0}%
_{5},L_{3}^{0}\right) $ and ${\bf M}_{{\cal N}}^{0}{\bf (U)=}\left(
-x_{3}(\partial _{\gamma \delta }^{2}u_{3})_{\gamma ,\delta =1,2},{\bf 0}%
_{6}\right) $, where $L_{3}^{0}=\varphi _{c}$ for Dirichlet conditions, then
(\ref{521}) may be rewritten as%
\[
\left\{
\begin{array}{l}
\forall {\bf V}^{1}{\bf =(v}^{1},{\psi }^{1})\in {\bf W}^{1},\medskip \\
\dint\nolimits_{Z}{\bf M}^{1}{\bf (V}^{1}){\cal R}\text{ }{^{t}}{\bf M}^{1}%
{\bf (U}^{1})\text{ d}z=-\dint\nolimits_{Z}{\bf M}^{1}({\bf V}^{1}){\cal R}%
\text{ }^{t}\left( {\bf M}_{{\cal M}}^{0}{\bf (U)+M}_{{\cal N}}^{0}{\bf (U)}%
\right) \text{d}z\text{.}%
\end{array}%
\right.
\]%
\smallskip \noindent This proves that
\[
{\bf U}^{1}={\bf U}_{{\cal M}}^{\lambda \mu }s_{\lambda \mu }({\bf \bar{u}})+%
{\bf U}_{{\cal N}}^{\lambda \mu }\partial _{\lambda \mu }^{2}u_{3}+{\bf U}%
^{3}L_{3}^{0}
\]%
where ${\bf U}_{{\cal M}}^{\lambda \mu }$, ${\bf U}_{{\cal N}}^{\lambda \mu
}\ $and ${\bf U}^{3}$ have been defined in (\ref{94})-(\ref{95})-(\ref{96}).
Moreover,%
\[
{\bf M}^{1}{\bf (U}^{1})={\bf M}^{1}{\bf (U}_{{\cal M}}^{\lambda \mu
})s_{\lambda \mu }(\overline{{\bf u}}{\bf )}+{\bf M}^{1}{\bf (U}_{{\cal N}%
}^{\lambda \mu })\partial _{\lambda \mu }^{2}u_{3}+{\bf M}^{1}{\bf (U}%
^{3})L_{3}^{0}\text{,}
\]%
which implies that%
\[
\begin{array}{l}
{\bf M}^{0}({\bf U})+{\bf M}^{1}{\bf (U}^{1})=\left( {\bf E}^{\lambda \mu }+%
{\bf M}^{1}{\bf (U}_{{\cal M}}^{\lambda \mu })\right) s_{\lambda \mu }(%
\overline{{\bf u}})\medskip \\
\ \ \ \ \ \ \ \ \ \ +\left( -x_{3}{\bf E}^{\lambda \mu }+{\bf M}^{1}{\bf (U}%
_{{\cal N}}^{\lambda \mu })\right) \partial _{\lambda \mu }^{2}u_{3}+\left( (%
{\bf 0}_{9},1)+{\bf M}^{1}{\bf (U}^{3})\right) L_{3}^{0}.%
\end{array}%
\]%
The choice in (\ref{93}) of ${\bf V}^{1}=(\widehat{{\bf v}}^{1},\psi
^{1})\in {\bf W}^{1}$ such that%
\[
{\bf M}^{1}{\bf (V}^{1})={\bf M}^{1}{\bf (u}_{{\cal M}}^{\alpha \beta
},\varphi _{{\cal M}}^{\alpha \beta })s_{\alpha \beta }(\overline{{\bf v}}%
{\bf )}+{\bf M}^{1}{\bf (u}_{{\cal N}}^{\alpha \beta },\varphi _{{\cal N}%
}^{\alpha \beta })\partial _{\alpha \beta }^{2}v_{3}
\]%
in the case of Dirichlet conditions and
\[
{\bf M}^{1}{\bf (V}^{1})={\bf M}^{1}{\bf (u}_{{\cal M}}^{\alpha \beta
},\varphi _{{\cal M}}^{\alpha \beta })s_{\alpha \beta }(\overline{{\bf v}}%
{\bf )}+{\bf M}^{1}{\bf (u}_{{\cal N}}^{\alpha \beta },\varphi _{{\cal N}%
}^{\alpha \beta })\partial _{\alpha \beta }^{2}v_{3}+{\bf M}^{1}{\bf (u}_{%
{\cal M}}^{3},\varphi _{{\cal M}}^{3})\widetilde{L}_{3}^{0}
\]%
in the case of mixed conditions, where ${\bf (v,}\widetilde{L}_{3})$ belongs
to ${\bf V}_{KL}\times \{0\}$ for Dirichlet conditions, to ${\bf V}%
_{KL}\times L^{2}(\omega )$ for local mixed conditions, and ${\bf V}%
_{KL}\times H^{1}(\omega )$ for nonlocal mixed conditions, completes the
proof. $\ \ \ \ \Box $

\pagebreak

\begin{center}
{\bf References}
\end{center}

\begin{enumerate}
\item Allaire G. Homogenization and two scale-convergence. {\it SIAM J. of
Math. An.} 1992; {\bf 23} (26); 1482-1518.

\item Bensoussan A, Lions J-L, Papanicolaou G. {\it Asymptotic Analysis for
Periodic Structures}. North Holland, Amsterdam; 1978.

\item Brassart M., Lenczner M. A two-scale model for the periodic
homogenization of the wave equation. {\it Journal de Math\'{e}matiques Pures
et Appliqu\'{e}es} 2010; {\bf 93(5)}, 474-517.

\item Caillerie D. Thin elastic and periodic plates. {\it Math. Meth. in
Applied Sciences} 1984; {\bf 6}; 159-191.

\item Canon \'{E}, Lenczner M. Models of elastic plates with piezoelectric
inclusions, part I: models without homogenization. {\em Math. Comp. Model.}
1997; {\bf 26} (5); 79-106.

\item Canon \'{E}, Lenczner M. Deux mod\`{e}les de plaques minces avec
inclusions pi\'{e}zo\'{e}lec-triques et circuits \'{e}lectroniques distribu%
\'{e}s. {\it C.R. Acad. Sci. Paris, S\'{e}rie II} 1998; {\bf 326}; 793-798.

\item Canon \'{E}, Lenczner M. Modelling of thin elastic plates with small
piezoelectric inclusions and distributed electronic circuits. Models for
inclusions that are small with respect to the thickness of the plate. {\it %
Journal of Elasticity} 2000; {\bf 55}; 111-141.

\item Casadei F., Ruzzene M., Dozio L., Cunefare K. A. Broadband vibration
control through periodic arrays of resonant shunts: experimental
investigation on plates. {\it Smart Materials and Structures} 2010; {\bf %
19(1)}, 015002.

\item Ciarlet P-G. {\it Mathematical Elasticity, Vol. II: Theory of Plates}.
North Holland, Amsterdam; 1997.

\item Ciarlet P-G, Destuynder P. A justification of the two-dimensional
plate model. {\it Comp. Meth. Appl. Mech. Eng.} 1979; {\bf 18}; 227-258.

\item Collet M., Ouisse M., Ichchou M., Ohayon R. Semi-active optimization
of 2D wave dispersion into shunted piezo-composite systems for controlling
acoustic interaction. {\it Smart Materials and Structures} 2012; {\bf 21(9)}%
; 094002.

\item Collet M., Ouisse M., Ichchou M. Structural energy flow optimization
through adaptive shunted piezoelectric metacomposites. {\it Journal of
Intelligent Material Systems and Structures} 2012; {\bf 23(15)}; 1661-1677.

\item Figueiredo I.N., Leal C. F. A piezoelectric anisotropic plate model.
{\em Asymptotic Analysis} 2005; {\bf 3-4 }; 327-346.

\item Ghergu M, Griso G, Mechkour H, Miara B. Homogenization of thin
piezoelectric perforated shells. {\em \ ESAIM, M2AN }2007; {\bf 41 }(5);
875-895.

\item Lenczner M. Homog\'{e}n\'{e}isation d'un circuit \'{e}lectrique. {\em %
C.R. Acad. Sci. Paris, s\'{e}rie II} 1997; {\bf 324(9)}, 537-542.

\item Lenczner M., Senouci-Bereksi, G. Homogenization of electrical networks
including voltage-to-voltage amplifiers. {\it Mathematical Models and
Methods in Applied Sciences} 1999; {\bf 9(06)}, 899-932.

\item Lenczner M., Mercier D. Homogenization of periodic electrical networks
including voltage to current amplifiers. {\it Multiscale Modeling \&
Simulation} 2004; {\bf 2(3)}, 359-397.

\item Licht C, Weller T. Mathematical Modeling of smart material and
structures. {\em East-West J. Math. special volume} 2007; 31-49.

\item Nguentseng G. A general convergence result for a functional related to
the theory of homogenization. {\it SIAM J. of Math. Anal.} 1989; {\bf 20}
(3); 608-623.

\item Raoult A. Contribution \`{a} l'\'{e}tude des mod\`{e}les d'\'{e}%
volution de plaques et \`{a} l'appro-ximation d'\'{e}quations d'\'{e}%
volution lin\'{e}aires du second ordre par des m\'{e}thodes multi-pas. Th%
\`{e}se de 3\`{e}me cycle, Universit\'{e} Pierre et Marie Curie, Paris, 1980.

\item Ratier N. Analog computing of partial differential equations. {\it %
IEEE 6th International Conference on in Sciences of Electronics,
Technologies of Information and Telecommunications} 2012; 275-282.

\item Sanchez-Palencia E. {\it Non-Homogeneous Media and Vibration Theory}.
Lecture Notes in Physic 127, Springer Verlag; 1980.

\item Sene A. Modelling of piezoelectric static thin plates. {\it Asymptotic
Analysis }2001; {\bf 127}; 1-20.

\item Simon J. Primitives de distributions et applications, {\em S\'{e}%
minaire d'Analyse}, Univ. Blaise Pascal 1991 6-7.

\item Thorp O., Ruzzene M., Baz, A. Attenuation and localization of wave
propagation in rods with periodic shunted piezoelectric patches. {\it Smart
Materials and Structures} 2001, 10(5), 979.

\item Weller T, Licht C. Analyse asymptotique de plaques minces lin\'{e}%
airement pi\'{e}zo-\'{e}lectriques. {\em C.R. Acad. Sci. Paris, s\'{e}rie I}
2002; {\bf 335}; 309-314.

\item Zemanian A-H. {\it Infinite Electrical Networks.} Cambridge University
Press, New York; 1991.

\item Zemanian A-H. Transfinite graphs and electrical networks. {\it Trans.
Amer. Math. Soc. }1992; {\bf 334}; 1-36.
\end{enumerate}

\end{document}